\documentclass[3p,12pt]{elsarticle}
\usepackage[english]{babel}
\usepackage{dsfont}
\usepackage{amssymb}
\usepackage{amsthm}
\usepackage{amsmath}
\usepackage{lineno}
\usepackage{empheq}
\usepackage{graphicx}
\usepackage{booktabs}
\usepackage{lineno}
\graphicspath{{picc/}}

\newdefinition{rmk}{Remark}
\newproof{pf}{Proof}
\newproof{pot}{Proof of Theorem \ref{thm2}}
\usepackage{caption}
\usepackage{subcaption}
\graphicspath{{picc/}}
\usepackage{amsmath}

\DeclareMathOperator*{\argmin}{arg\,min}
\begin{document}
\begin{frontmatter}
\title{Generating the Laguerre expansion coefficients by solving a one-dimensional transport equation.}
\author{Andrew V. Terekhov}
\ead{andrew.terekhov@mail.ru}
\address{Institute of Computational Mathematics and Mathematical Geophysics, 630090, Novosibirsk, Russia}
\address{Novosibirsk State Technical University, 630073, Novosibirsk, Russia}
\begin{abstract}
Spectral methods based on integral transforms may be efficiently used to solve differential equations in some special cases. This paper considers a different approach in which algorithms are proposed to calculate integral Laguerre transform by solving a one-dimensional transport equation. In contrast to the direct calculation of improper integrals of rapidly oscillating functions, these procedures make it possible to calculate the expansion coefficients of a Laguerre series expansion with better stability, higher accuracy, and less computational burden.
\end{abstract}
\end{frontmatter}
\section{Introduction}
The Laguerre integral transform has been used in various fields of mathematical simulation to solve acoustics and elasticity equations \cite{fatab2011,Terekhov:2013,Terekhov2015206,Mikhailenko1999}, Maxwell and heat conduction equations \cite{Mikhailenko2008,Colton1984}, and spectroscopy problems \cite{Jo2006}. The Laguerre transform has proved to be a very efficient tool in constructing a stable algorithm of wave field continuation when solving inverse problems of seismic prospecting \cite{Terekhov2017,Terekhov2018} and many others. The Laguerre transform has served as a basis for the development of numerical methods of inversion of Laplace \cite{Weeks1966,Abate1996,Strain1992} and Fourier \cite{Weber1980} transforms. In numerically solving differential equations by applying the Laguerre transform in time and approximating space derivatives one has to solve definite well-conditioned systems of linear algebraic equations. For the latter one can use fast convergent algorithms of computational linear algebra \cite{Golub1989,Samarski_Nikolaev}. In addition, in contrast to the Fourier transform, to calculate the coefficients of the Laguerre series one and the same operator, which does not depend on the number of the harmonic being calculated, is inverted several times. On the contrary, the operator obtained by the Fourier transform will depend on the frequency. This property of the Laguerre transform allows using efficient parallel preconditioning procedures to solve systems of linear algebraic equations, for instance, on the basis of the dichotomy algorithm \cite{Terekhov:2013,terekhov:Dichotomy,Terekhov2016}, which was specially developed to invert one and the same matrix for different right-hand sides.

Consider Laguerre functions \cite{abramowitz+stegun}, which are defined as
\begin{equation}
l_n(t)=e^{-t/2}L_n(t),\quad t\geq 0
\label{laguerre_function}
\end{equation}
where $L_n(t)$ is the Laguerre polynomial of degree $n$, which is defined by the Rodrigues formula
$$
L_n(t)=\frac{ e^{t}}{n!}\frac{d^n}{dt^n}\left(t^ne^{-t}\right)=\frac{1}{n!}\left(\frac{d}{dt}-1\right)^nt^n.
$$
We will use $L_2[0,\infty)$ to denote the space of square integrable functions $f:[0,\infty)\rightarrow \mathbb{R}$
$$
L_2[0,\infty)=\left\{f:\int_{0}^{\infty}|f(t)|^2dt< \infty\right\}.
$$
The Laguerre functions are a complete orthonormal system in $ L_2[0,\infty)$
\begin{equation}
  \int_{0}^{\infty}l_m(t)l_n(t)dt=
  \left\{\begin{array}{ll}
  0,& m \neq n ,\\
  1,& m=n,
 \end{array}\right.
 \label{lag_orho}
 \end{equation}
This guarantees that for any function $f(t)\in L_2[0,\infty)$ there is a Laguerre expansion
\begin{subequations}
\label{series_lag}
\begin{empheq}{align}
\label{series_lag.sum}
f(t)\sim\eta\sum_{m=0}^{\infty}\bar{a}_m l_m(\eta t), \quad t  \geq 0,\  \eta > 0,\\ \bar{a}_m=\int_{0}^{\infty}f(t)l_m(\eta t)dt ,
\label{series_lag.int}
\end{empheq}
\end{subequations}
where $\eta$ is a scaling parameter for the Laguerre functions to increase the convergence rate of the series (\ref{series_lag.sum}).

The Laguerre function values for large $n$ are bounded from above, since an asymptotic representation \cite{Szegoe1975}
\begin{equation}
\label{eq:assymptotic}
l_n(t)=\frac{1}{\pi^{1/2}(nt)^{1/4}}\left(\cos(2\sqrt{nt}-\pi/4)\right)+O\left(\frac{1}{n^{3/4}}\right),\quad t\in[a,b],\quad 0< a< b < \infty,
\end{equation}
is valid. However, one of the problems of numerical implementation of the transform  (\ref{series_lag.int}) is that in calculating the Laguerre functions for  $t>1$ the values of the Laguerre polynomials $L_n(t)$ rapidly increase with increasing $n$,  which leads to an error of "overflow"{}. On the contrary, in calculating the multiplier $\exp({-t/2})$ there may be an error of "underflow"{}. For small $n$ and $t$ the Laguerre functions can be calculated by the formula
\begin{equation}
l_n(t)=\left[e^{-t/4}\tilde{L}_n(t)\right]e^{-t/4}.
\label{eq:little_laguerre_function}
\end{equation}
Specifically, first we calculate the expression in the square brackets using a second order recurrence formula \cite{Rainville1971}:
\begin{equation}
\begin{array}{l}
(n+1)\tilde{L}_{n+1}(t)=(2n+1-t)\tilde{L}_n(t)-n\tilde{L}_{n-1}(t), \quad n\geq 1,\\\\
\tilde{L}_1(t)=(1-t)e^{-t/4},\quad \tilde{ L}_0(t)=e^{-t/4}.
\end{array}
\label{recurrence_laguerre}
\end{equation}
Then, multiplying the result by the second exponential multiplier, we calculate the Laguerre function. If the calculations are made at $128$-bit real computer precision, this method excludes situations of the "overflow"{} and "underflow"{} types for $n$-values that do not exceed several thousand and $t<20$, $\eta<1200$. Since high-precision arithmetic is used, as a rule, with software (but not hardware), the use of high precision considerably decreases the efficiency of calculations. Therefore, to save the calculation time $128$-bit arithmetic should only be used to calculate the Laguerre functions, whereas the summation in approximating the integral (\ref{series_lag.int}) can be made using standard $64$-bit precision.

Another problem of implementing the Laguerre transform is caused by the fact that the Laguerre functions of the $n$-th order on the interval $ 0<t<4n$ oscillate \cite{Temme1990}, and the strongest oscillations are near zero (see Fig.~\ref{pic:laguerre_function}). This brings up a problem of finding a method to integrate rapidly oscillating functions. To overcome this difficulty, an algorithm to calculate the integral (\ref{series_lag.int}) is proposed in \cite{Litko1989}. This algorithm is based on quadratures of high-order accuracy, which make it possible to calculate the Laguerre series expansion coefficients whose number $n$ is not greater than several hundreds. However, one should take into account that the quadratures of high orders are defined on nonuniform grids, which may not allow their use if the function to be approximated is given in the form of a time series for equal-spaced time intervals. For analytical functions, approaches based on the Laplace transform and Cauchy's integral formula can be used \cite{Abate1996}:
\begin{figure}[!htb]
\centering
\begin{subfigure}[b]{0.47\textwidth}
\includegraphics[width=\textwidth]{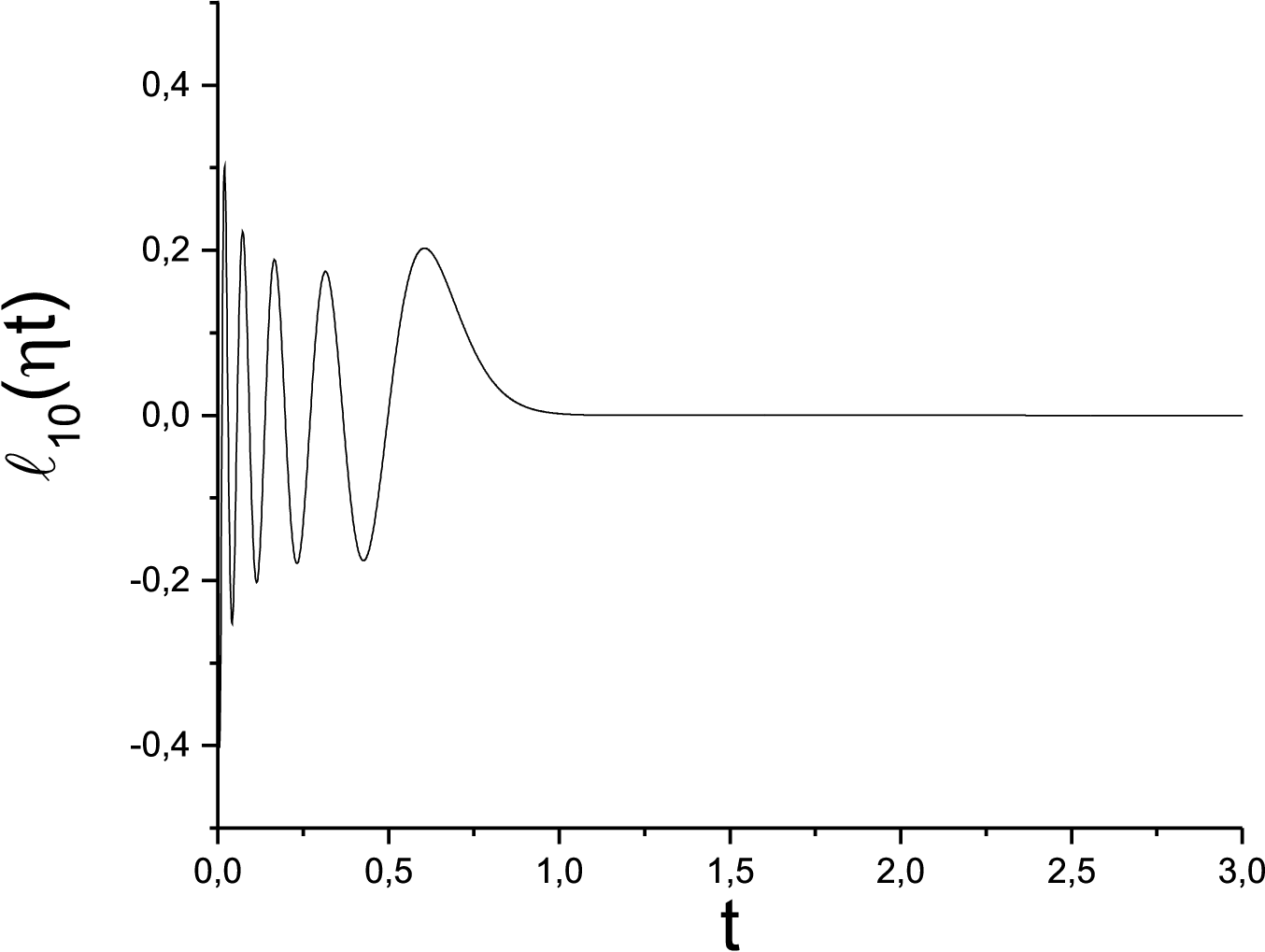}
\end{subfigure}
~
\begin{subfigure}[b]{0.47\textwidth}
\includegraphics[width=\textwidth]{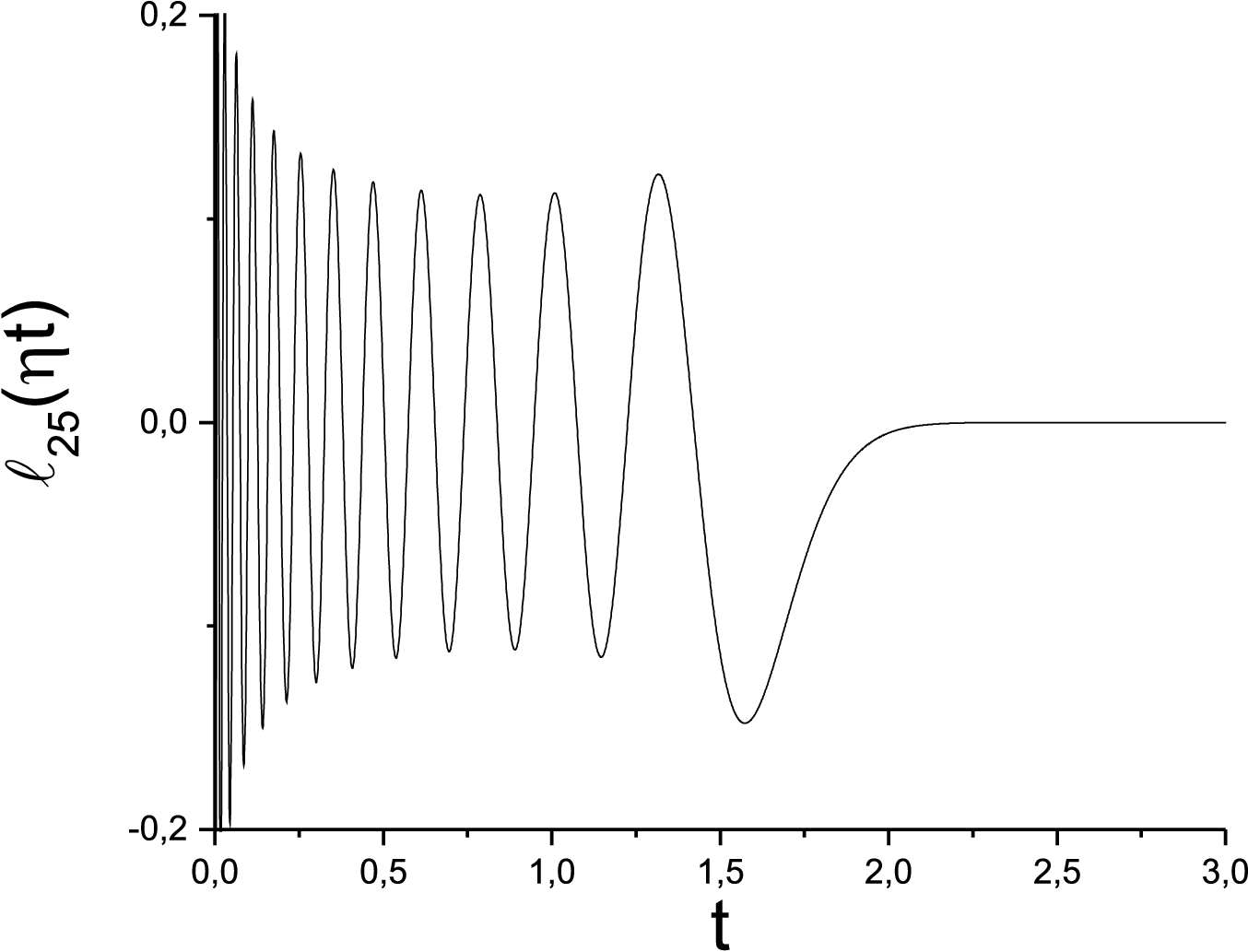}
\end{subfigure}
\begin{subfigure}[b]{0.47\textwidth}
\includegraphics[width=\textwidth]{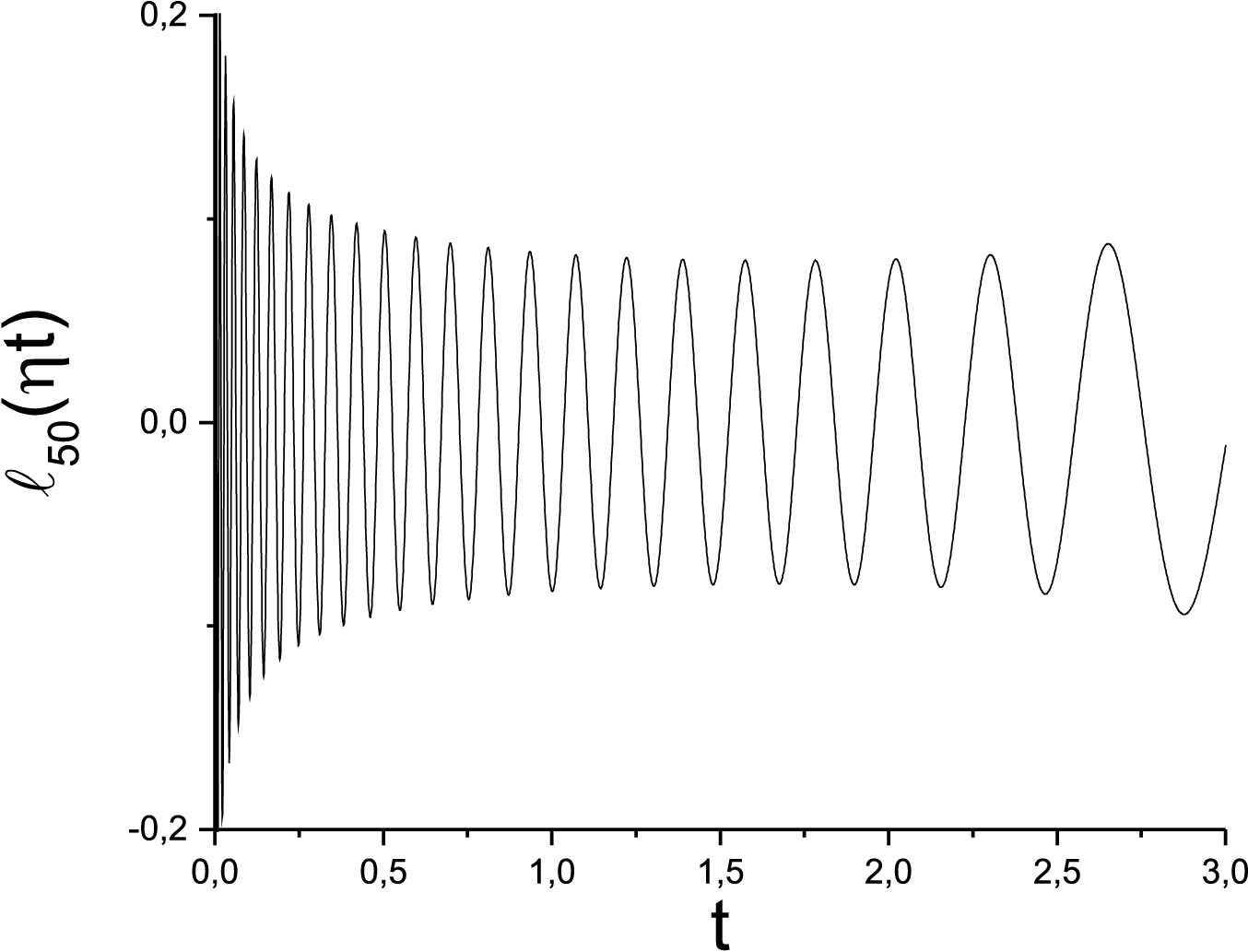}
\end{subfigure}
~
\begin{subfigure}[b]{0.47\textwidth}
\includegraphics[width=\textwidth]{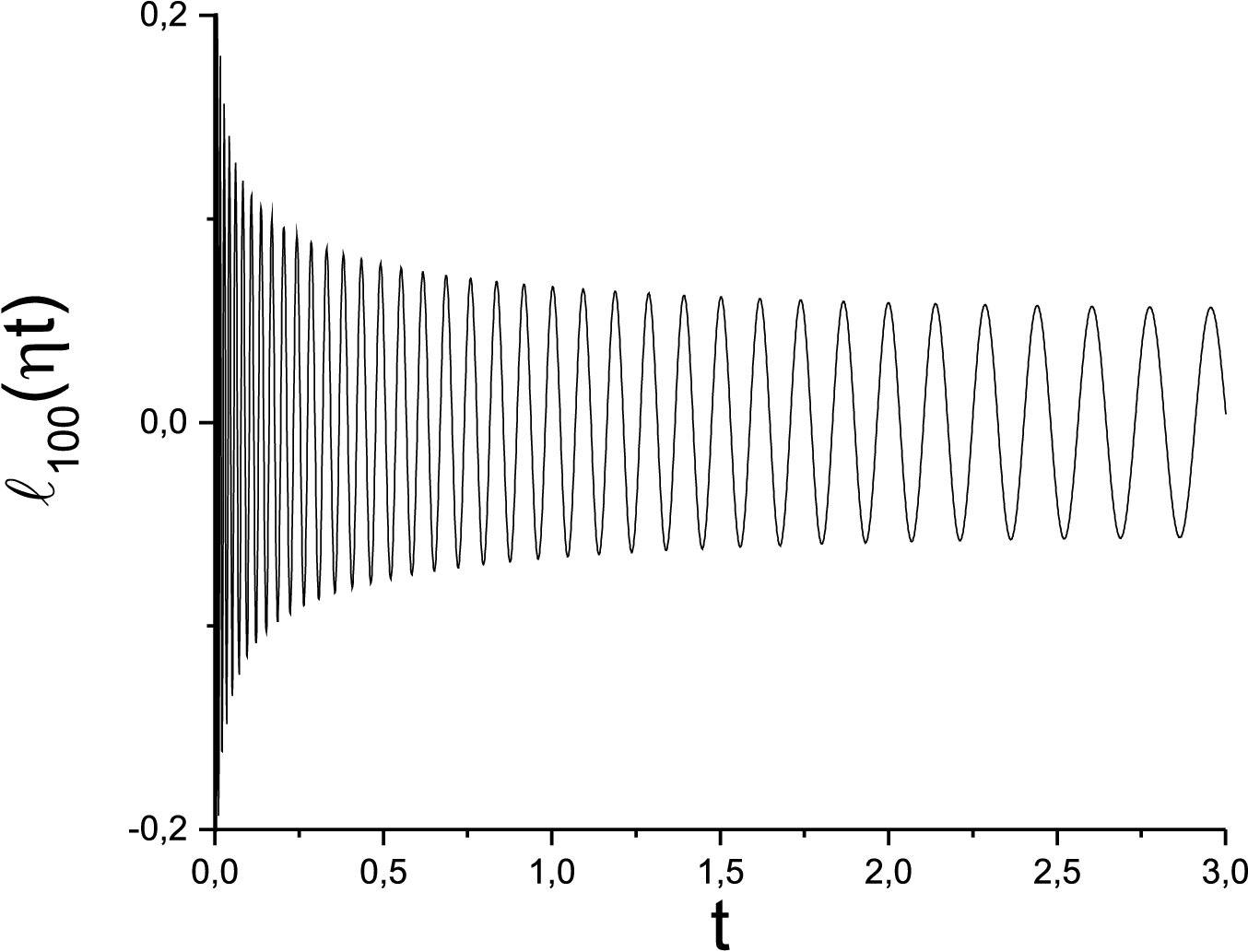}
\end{subfigure}
\caption{Laguerre functions of various orders for the transformation parameter \mbox{$\eta=60$}.}
\label{pic:laguerre_function}
\end{figure}
\begin{equation}\label{lag_lap}
\bar{a}_n=\frac{1}{2\pi \mathrm{I}}\int_{C_r}\left[\frac{\hat{f}((1+z)/2(1-z))}{1-z}\right]z^{-(n+1)}dz.
\end{equation}
Here the expression in the square brackets is a generating Laguerre function, which is analytical in the circle $C_r$ of radius $r$ with the center at the origin of coordinates, and $\hat{f}(s)$ is the Laplace transform for the function $f(t)$. The imaginary unit is denoted by $\mathrm{I}=\sqrt{-1}$. In calculating the integral (\ref{lag_lap}), the necessary preliminary Laplace transform for a function expanded into a Laguerre series makes it difficult to use this algorithm.

Another method of calculating the expansion coefficients is given by an integral of the form \cite{Weber1980}
\begin{equation}
\label{lag_four}
\bar{a}_n=\frac{1}{2\pi}\int_{0}^{2\pi}\left[\frac{1}{2}\left(1+\mathrm{I}\cot\frac{z}{2}\right)f\left(\frac{1}{2}\cot\frac{z}{2}\right)\right]e^{-\mathrm{i}nz}dz.
\end{equation}
Note that the cotangent function has singularities at points $0$ and $\pi$. This complicates the calculation of the expansion coefficients if the function being approximated is discrete. One more method based on the Laplace and Fourier transforms was proposed in \cite{Weeks1966}. Thus, the above-mentioned approaches are most suited for approximating smooth analytical functions, for which the Laplace or Fourier transform is known. However, in solving applied problems the initial data may be specified in the form of time series with low smoothness, which calls for the development of additional procedures for this case.

In this paper, a new method to calculate the Laguerre series coefficients is proposed. It is based on solving a one-dimensional transport equation. This is a distinguishing feature of the approach, since the integral transforms are used, as a rule, to solve differential equations. In contrast to this, the one-dimensional transport equation is solved to implement the integral Laguerre transform. With this approach to the problem, a stable and rather accurate algorithm which is less expensive than the direct calculation of the integral (\ref{series_lag.int}) is proposed. In addition, an efficient variant of the method will be considered for the approximation of functions on large intervals.

\section{Expansion algorithms}
\subsection{Main formulas}
Consider the following initial boundary value problem for a one-dimensional transport equation:
\begin{equation}
\left\{
\begin{array}{ll}
\displaystyle \frac{\partial v}{\partial t }-\frac{\partial v}{\partial x }=0, \quad t>0,\quad  -\infty< x < +\infty,\\\\
\displaystyle v(x,0)=f(x).
\end{array}\right.
\label{advection_eq}
\end{equation}
On taking the Laguerre transform in time of the problem (\ref{advection_eq}), it can be written in the form \cite{Mikhailenko1999}
\begin{equation}
\label{advection_laguerre}
\left(\frac{\eta}{2}-\partial_x\right)\bar{v}_m=-\Phi(\bar{v}_m),
\end{equation}
where
\begin{equation}
\Phi(\bar{v}_m)=-f+\eta\sum_{j=0}^{m-1}\bar{v}_j.
\label{phi_function}
\end{equation}
Taking into consideration $$\Phi(\bar{v}_{m})=\eta \bar{v}_{m-1}+\Phi(\bar{v}_{m-1}),$$
let us turn to another form of (\ref{advection_laguerre})
\begin{subequations}
\label{advection_laguerre3}
  \begin{empheq}[left=\empheqlbrace]{align}
  \label{advection_laguerre3.a}
  \left(\frac{\eta}{2}-\partial_x\right)\bar{v}_{0}-f=0,\\
  \label{advection_laguerre3.b}
  \left(\frac{\eta}{2}-\partial_x\right)\bar{v}_{m}=\left(-\frac{\eta}{2}-\partial_x\right)\bar{v}_{m-1},\ m=1,2,...
  \end{empheq}
\end{subequations}
Then, taking the Fourier transform in the variable $x$, we have
\begin{subequations}
\label{advection_laguerre2}
  \begin{empheq}[left=\empheqlbrace]{align}
  \label{advection_laguerre2.a}
  \left(\frac{\eta}{2}-\mathrm{I}k\right)\bar{V}_{0}(k)-\tilde{f}(k)=0,\\
  \label{advection_laguerre2.b}
  \left(\frac{\eta}{2}-\mathrm{I}k\right)\bar{V}_{m}(k)=\left(-\frac{\eta}{2}-\mathrm{I}k\right)\bar{V}_{m-1}(k),\ m=1,2,...,
  \end{empheq}
\end{subequations}
where $k$ is the wavenumber. Expressing the sought-for function in explicit form, we have
\begin{equation}
\begin{array}{c}\displaystyle
\bar{V}_{m}(k)=\tilde{f}(k){\left({-\frac{\eta}{2}-\mathrm{I}k}\right)^m}/{\left({\frac{\eta}{2}-\mathrm{I}k}\right)^{m+1}}.
\end{array}
\label{lag_fourier}
\end{equation}

Again, consider the problem (\ref{advection_eq}), but with periodic boundary conditions of the form \mbox{$v(0,t)=v(T,t)$}, where $T$ determines the boundary of the interval of approximation of the function $f(t)$, $t \in [0,T]$. In this case the solution to equation (\ref{advection_laguerre})  has the form of summation of solutions of the form (\ref{lag_fourier}) for a discrete set of frequencies, $k_j=2 \pi j/T$, $j=0,1,...,N_x$:
\begin{equation}
\label{main_formula}
\bar{v}_{m}(p)\approx\sum_{j=0}^{N_x}\tilde{V}_{m}(k_{j})\exp\left(\mathrm{I}\frac{2\pi j p}{T}\right).
\end{equation}
Subject to the solution (\ref{main_formula}) for the transport equation, the function $f(t)$, given as an initial condition, will move in the direction $x=0$. By writing the solution to the transport equation at the point $x=0$, we see that the sought-for coefficients of the expansion (\ref{series_lag.int}) for the function $f(t)$ can be calculated as $\bar{a}_m=\bar{v}_{m}(0)$.

Although the expansion coefficients are calculated by formulas (\ref{lag_fourier}), (\ref{main_formula}) with $O(nN_x)$ arithmetic operations, that is, the algorithm is not fast, the above method, proposed for implementing the Laguerre transform, has some important advantages over the direct calculation of the improper integral of the rapidly oscillating function (\ref{series_lag.int}). First, from the definition of the absolute value of a complex number we have the identity
$$\left|{\left({-\frac{\eta}{2}-\mathrm{I}k}\right)}/{\left({\frac{\eta}{2}-\mathrm{I}k}\right)^{}}\right|\equiv1,$$
which guarantees stability of the calculation and the absence of "overflow"{} or "underflow"{} situations for any $T$ and $n$, which is a problem in calculating the Laguerre functions by formula (\ref{laguerre_function}). Second, as it will be shown below, the calculations by formulas (\ref{lag_fourier}),(\ref{main_formula}) can be made with single $32$-bit real precision, which increases the accuracy of the calculations by using higher vectorization. On the contrary, the considerable spread in the Laguerre function values calls for $64$-bit precision calculations. Third, despite the presence of strong oscillations of the Laguerre functions at the origin of coordinates (Fig.~\ref{pic:laguerre_function}), the spectral approach does not require using nonuniform grids or quadratures of high-order accuracy to retain a given accuracy on the entire approximation interval. From a practical viewpoint, it is much more convenient to specify the number of harmonics $N_x$ of the Fourier series instead of the grid size, since the boundaries of the spectrum of the function being approximated are, as a rule, either known beforehand or can be determined in an efficient way.

A shortcoming of the computational model being considered is that this method of calculating the expansion coefficients of the Laguerre series adds a fictitious periodicity of the form $f(t)=f(t+bT)$, where b is any nonnegative integer. To remove the undesirable periodicity, two fundamentally different approaches will be proposed below.
\subsection{Energy-dependent truncation of Laguerre Series}
\label{section:enetgy_truncation}
\begin{figure}[!htb]
\centering
\begin{subfigure}[b]{0.47\textwidth}
\includegraphics[width=\textwidth]{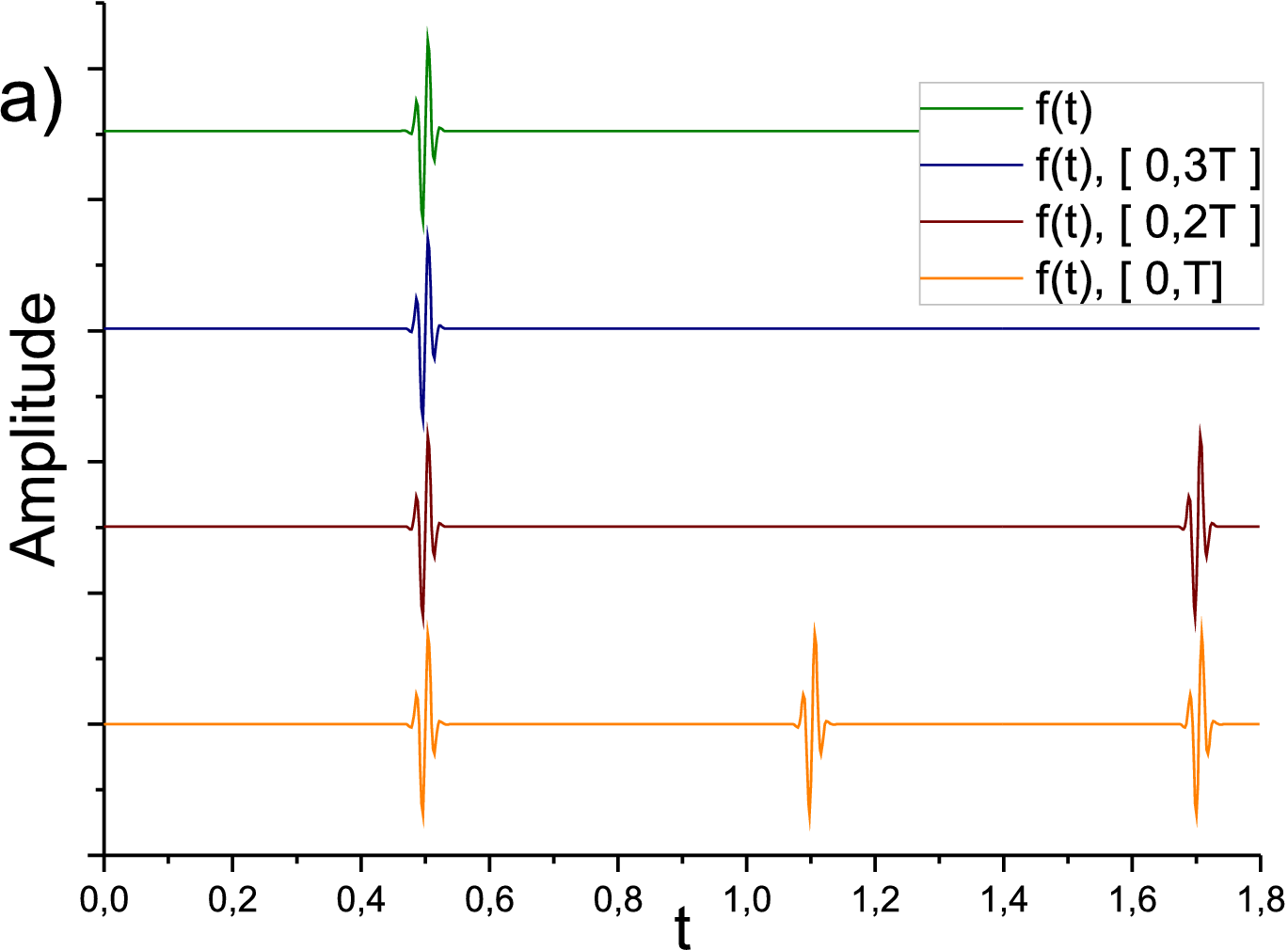}
\end{subfigure}
~
\begin{subfigure}[b]{0.47\textwidth}
\includegraphics[width=\textwidth]{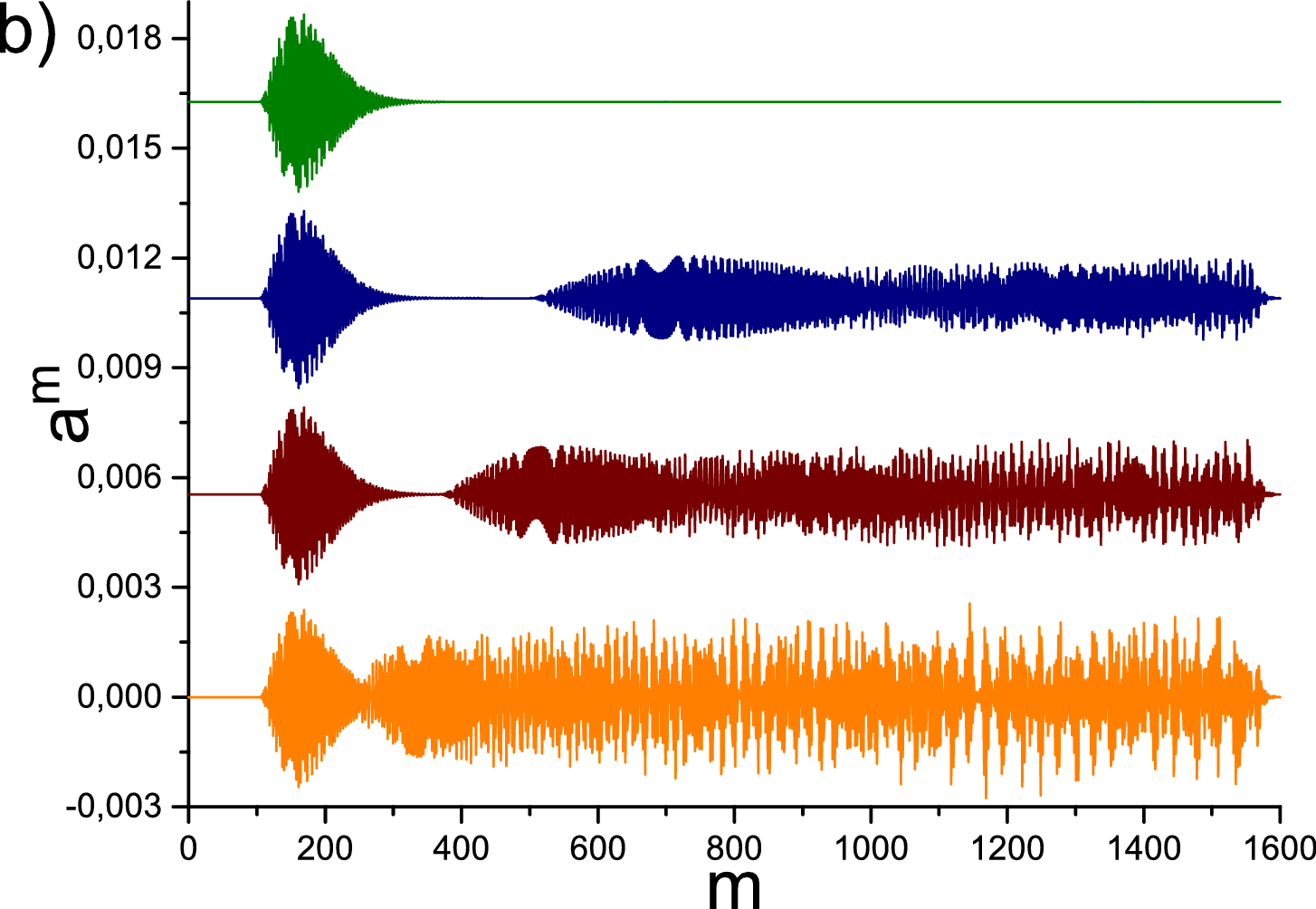}
\end{subfigure}
\caption{a) Function (\ref{source}) and its approximation by formula (\ref{main_formula}) for intervals of various lengths, b) Laguerre spectrum.}
\label{pic:conjg222}
\end{figure}

Consider an approach which allows removing the fictitious periodicity in the calculation of the expansion coefficients of the Laguerre series by formula (\ref{main_formula}). Fig.~\ref{pic:conjg222}b shows the Laguerre series expansion coefficients for the function $f(t)$ in Fig.~\ref{pic:conjg222}b specified by the formula
\begin{equation}
f(t)=\exp\left[-\frac{(2\pi f_0(t-t_0))^2}{g^2}\right]\sin(2\pi f_0(t-t_0)),
\label{source}
\end{equation}
where $t_0=0.5,\, g=4,\,f_0=30$. It is evident from Fig.~\ref{pic:conjg222}b that to exclude the undesirable periodicity it is sufficient to increase the calculation interval from $[0,T]$ to $[0,3T]$ assuming that for $t \in [T,3T]$ the function is zero. Then, once the expansion coefficients have been calculated, remove the coefficients with numbers $m>m_0 \approx 400$. In the calculations for smaller approximation intervals, for instance $[0,T]$ or $[0,2T]$, the fictitious periods of the function in the spectral domain cannot be separated, since the abrupt truncation of the series will cause oscillations in the entire approximation interval.

To automatically determine the number of the remaining expansion coefficients, we use Parseval's  relation
\begin{equation}
\int_{0}^{\infty} v^2(t) dt=\sum_{m=0}^{\infty}\left(\bar{v}_{m}\right)^2.
\label{pseudo_energy}
\end{equation}
On the basis of this relation the maximum number of the Laguerre series coefficients $m_0$ is determined from the condition
\begin{equation}
\argmin_{m_0}\left|{\int_{0}^{T} v^2(t) dt-\sum_{m=0}^{m_0}\left(\bar{v}_{m}\right)^2}{}\right|.
\label{pseudo_energy2}
\end{equation}
Now let us formulate an algorithm of expanding the function $f(t)$ in a Laguerre series.
\\\\
{\bf Algorithm 1} to approximate a function $f(t)$  on the interval $t\in\left[0,T\right]$ by a Laguerre series:
\begin{enumerate}
  \item  	Calculate $\tilde{f}=FFT(f)$ on the basis of a fast algorithm of the discrete Fourier transform.
  \item  	Calculate the expansion coefficients of the series (\ref{series_lag.sum}) by formula (\ref{main_formula}) and the equality $\bar{a}_m=\bar{v}_{m}(0)$.
  \item  On the basis of formula (\ref{pseudo_energy2}), leave intact only the first $m_0$ coefficients of the series (\ref{series_lag.sum}).
\end{enumerate}

The above-considered a posteriori method of removing the fictitious periodicity is not convenient from a practical viewpoint, since the spectra of nonsmooth functions may be rather large. This may not allow separating the first period of the function being approximated from the subsequent fictitious periods in the spectral domain. Also, in approximating functions of various smoothness it is not clear how many times the approximation interval must be increased to reliably remove the fictitious periodicity. In this case too great increase in the approximation interval length may cause a considerable increase in the computational costs. To solve these problems, an alternative procedure of removing the fictitious periodicity not requiring a posteriori analysis of the Laguerre spectrum will be developed.
\subsection{Shift and Conjugation procedures}
Let us develop two auxiliary procedures to modify the Laguerre series coefficients, which we call shift and conjugation. These will allow us to propose an alternative procedure of removing the fictitious periodicity, as well as a procedure of reducing the computational costs when a function is expanded in a series for large approximation intervals.

Consider an analytical solution to the following initial boundary value problem
\begin{equation}
\label{boundary3}
\left\{
\begin{array}{ll}
\displaystyle \frac{\partial v}{\partial t }+\frac{\partial v}{\partial x }=0,&   t>0,\  x>0,\\\\
v(0,t)=f(t), & t\geq 0,\\
v(x,0)=0, & x\geq 0,\\
f(0)=0.
\end{array}\right.
\end{equation}
As in solving the problem (\ref{advection_eq}), we again apply the Laguerre transform in time to the transport equation, and obtain the equation
\begin{equation*}
\label{advection_laguerre222}
\left(\frac{\eta}{2}+\partial_x\right)\bar{v}_m=-\Phi(\bar{v}_m).
\end{equation*}
For the boundary conditions (\ref{boundary3}) to be satisfied, we use the Laguerre transform but not the Fourier transform to calculate the functions $\bar{v}_m(x)$, that is, search for a solution of the form
\begin{equation}
\bar{v}_m(x)=\kappa \sum_{j=0}^{\infty}W_{m,j}l_j(\kappa x),\quad m=0,1,2...,
\label{laguerre_x_expansion}
\end{equation}
where the transformation parameter $\kappa>0$. Then, on applying the Laguerre spatial transform to equation  (\ref{boundary3}), we have
\begin{subequations}
\label{lagx}
\begin{empheq}[left=\empheqlbrace]{align}
\label{lagx.1}
\displaystyle \left(\eta+\kappa\right)W_{m,0}=\left(-\eta+\kappa\right)W_{m-1,0}+2\left(\bar{f}_m- \bar{f}_{m-1}\right),&\quad m=0,1,...,\\
\label{lagx.2}
\left(\eta+\kappa\right)W_{m,j}+2\Upsilon(W_{m,j})=\left(-\eta+\kappa\right)W_{m-1,j}+2\Upsilon(W_{m-1,j}),&\quad m=0,1,...;\, j=1,2,..,
\end{empheq}
\end{subequations}
where
\begin{equation}
\Upsilon\left(W_{m,j}\right)=\kappa \sum_{i=0}^{j-1}{W_{m,i}}=\kappa W_{m,j-1}+\Upsilon\left(W_{m,j-1}\right),
\label{eq1}
\end{equation}
 and
$W_{m,j}\equiv0,\;\bar{f}_{m}\equiv0,\quad \forall \; m<0.$
\\\\
Taking (\ref{eq1}) into account, equation (\ref{lagx.2})  takes the following form:
\begin{equation}
\left(\eta+\kappa\right)W_{m,j}+\left(\eta-\kappa\right){W}_{m-1,j}=
\left(\eta-\kappa\right)W_{m,j-1}+\left(\eta+\kappa\right)W_{m-1,j-1},\quad  m=0,1,...;\,j=1,2,...\ .
\end{equation}
Taking $\kappa=\eta$, we finally obtain
\begin{equation}
\label{lagx3}
\left\{\begin{array}{ll}
W_{m,0}=\left(\bar{f}_m-\bar{f}_{m-1}\right),& m=0,1,...,
\\
W_{m,j}=W_{m-1,j-1},& m=0,1,...;\; j=1,2,...
\end{array}
\right.
\end{equation}
Based on (\ref{series_lag.sum}), (\ref{laguerre_x_expansion}) and (\ref{lagx3}), the final solution to problem (\ref{boundary3})  in the time domain is as follows:
\begin{equation}
\label{solution1}
v(x,t)=\kappa\sum_{m=0}^{\infty} \left(\eta\sum_{j=0}^mW_{m-j,0}l_{j}(\kappa x)\right)l_{m}(\eta t).
\end{equation}
Changing the order of summation, we can also write
\begin{equation}
\label{solution2}
v(x,t)=\eta\sum_{j=0}^\infty\left(\kappa\sum_{m=0}^{\infty}W_{m,0}l_{m+j}(\eta t)\right)l_{j}(\kappa x).
\end{equation}
It follows from formulas (\ref{solution1}) and (\ref{solution2}) that the expressions in the brackets are the Laguerre series coefficients. Then, taking into account the relations (\ref{lagx3}), we introduce two transforms with a parameter $\tau\geq 0$:
\begin{subequations}
\begin{empheq}{align}
\label{shift_proc}
\mathbb{S}\left\{\bar{a}_m;\tau\right\}=\sum_{j=0}^m\left(\bar{a}_{m-j}-\bar{a}_{m-j-1}\right)l_{j}(\eta \tau)\ ,
\\ \mathbb{Q}\left\{\bar{a}_j;\tau\right\}=\sum_{m=0}^{\infty}\left(\bar{a}_m-\bar{a}_{m-1}\right)l_{m+j}(\eta \tau),\ \mathrm{where}\ \bar{a}_{-1}\equiv0.
\label{reverse_proc}
\end{empheq}
\label{transform2}
\end{subequations}

\begin{figure}[!htb]
\centering
\begin{subfigure}[b]{0.47\textwidth}
\includegraphics[width=\textwidth]{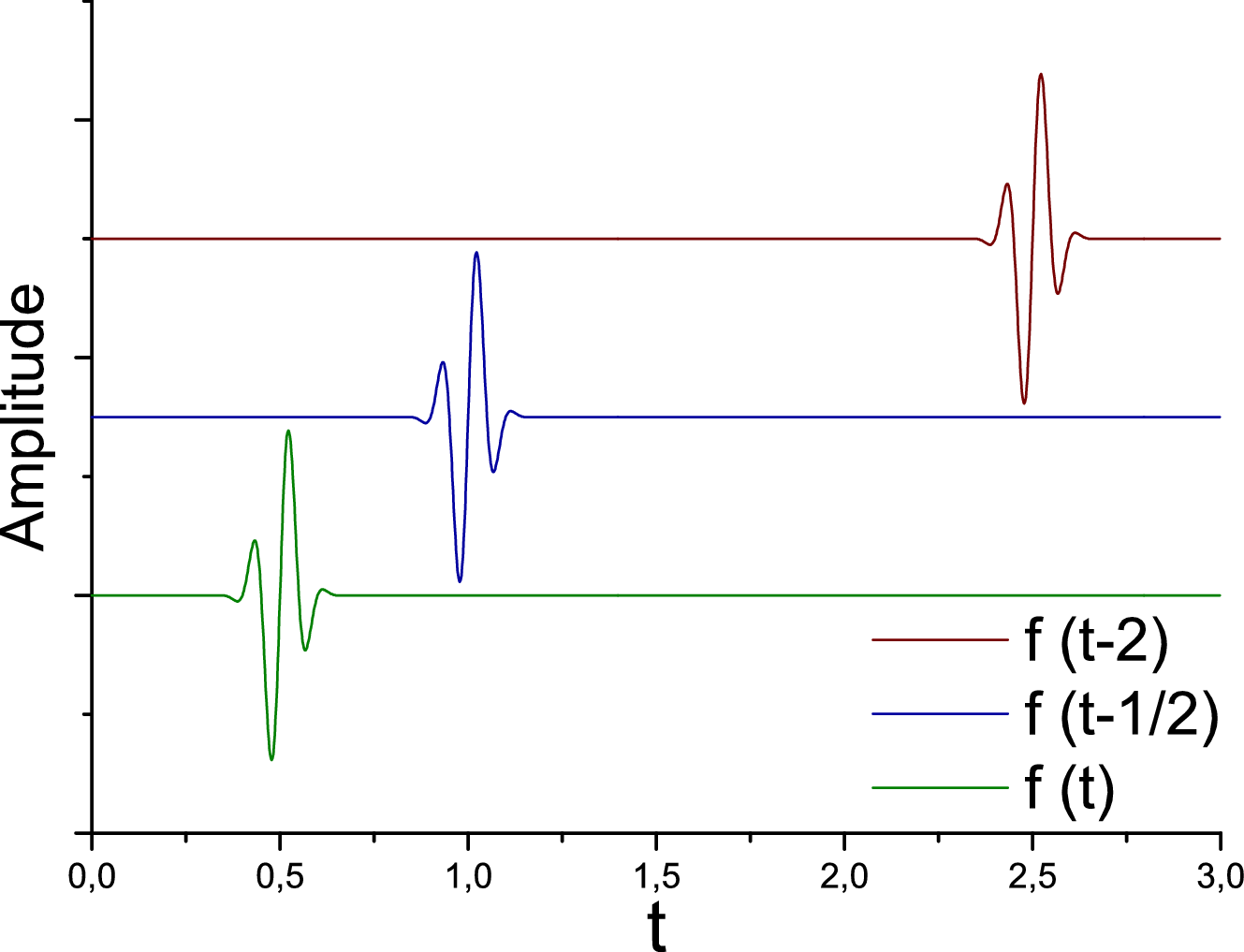}
\end{subfigure}
~
\begin{subfigure}[b]{0.47\textwidth}
\includegraphics[width=\textwidth]{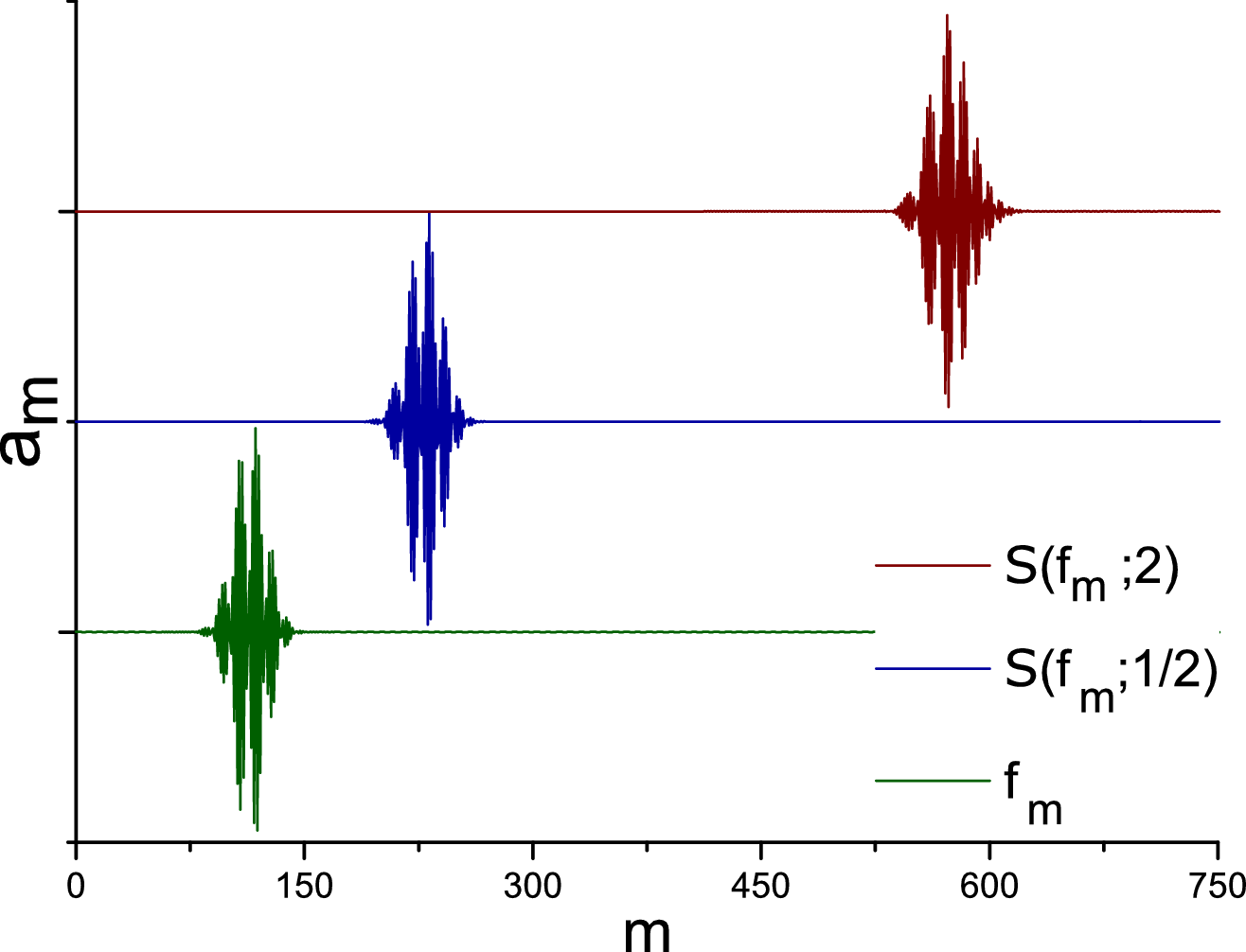}
\end{subfigure}
\caption{a) Function (\ref{source}) and b) Laguerre spectrum for various values of parameter $\tau$ of shift operator $\mathbb{S}\{\bar{f}_m;\tau\}$.}
\label{pic:shift_test}
\end{figure}

\begin{figure}[!htb]
\centering
\begin{subfigure}[b]{0.47\textwidth}
\includegraphics[width=\textwidth]{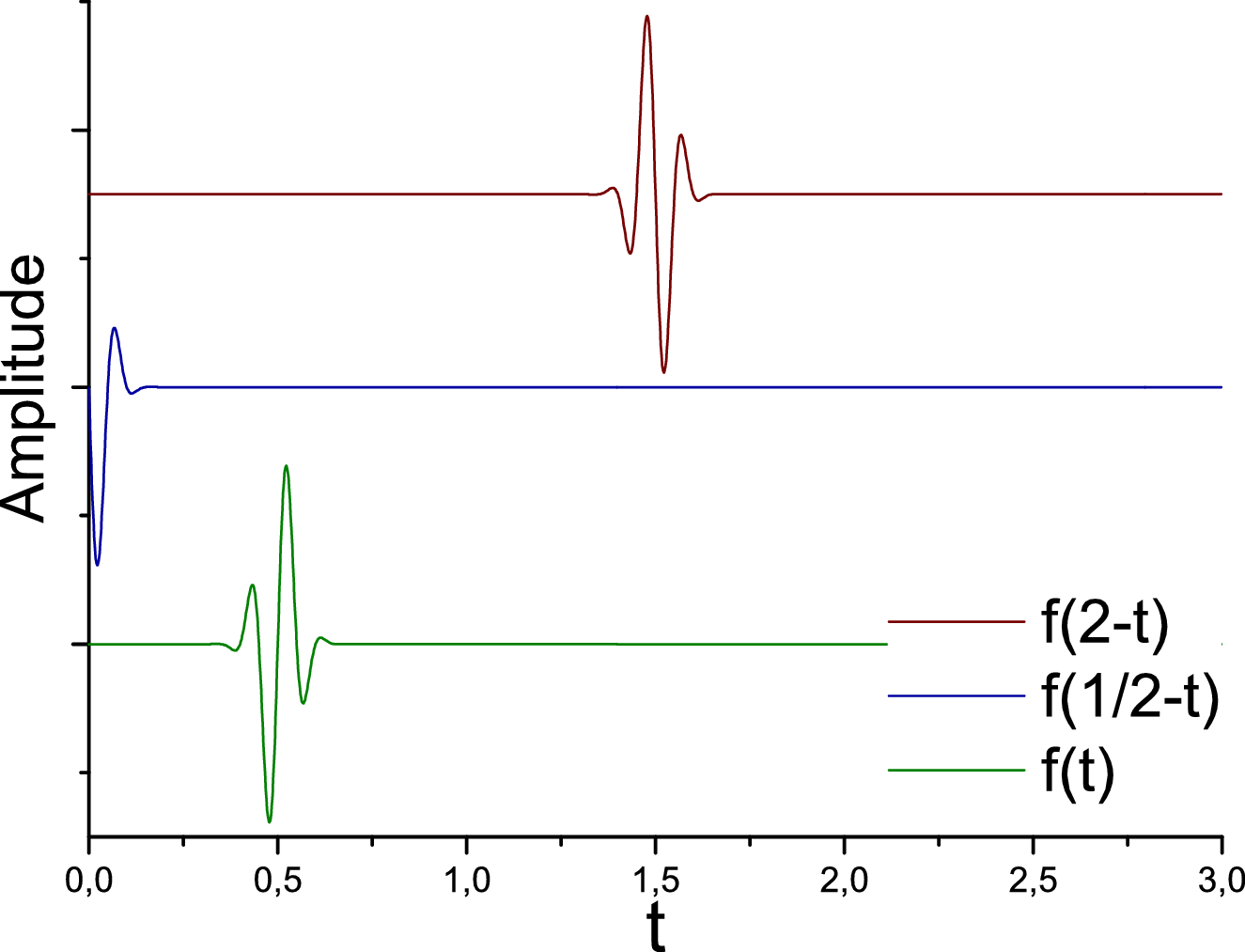}
\end{subfigure}
~
\begin{subfigure}[b]{0.47\textwidth}
\includegraphics[width=\textwidth]{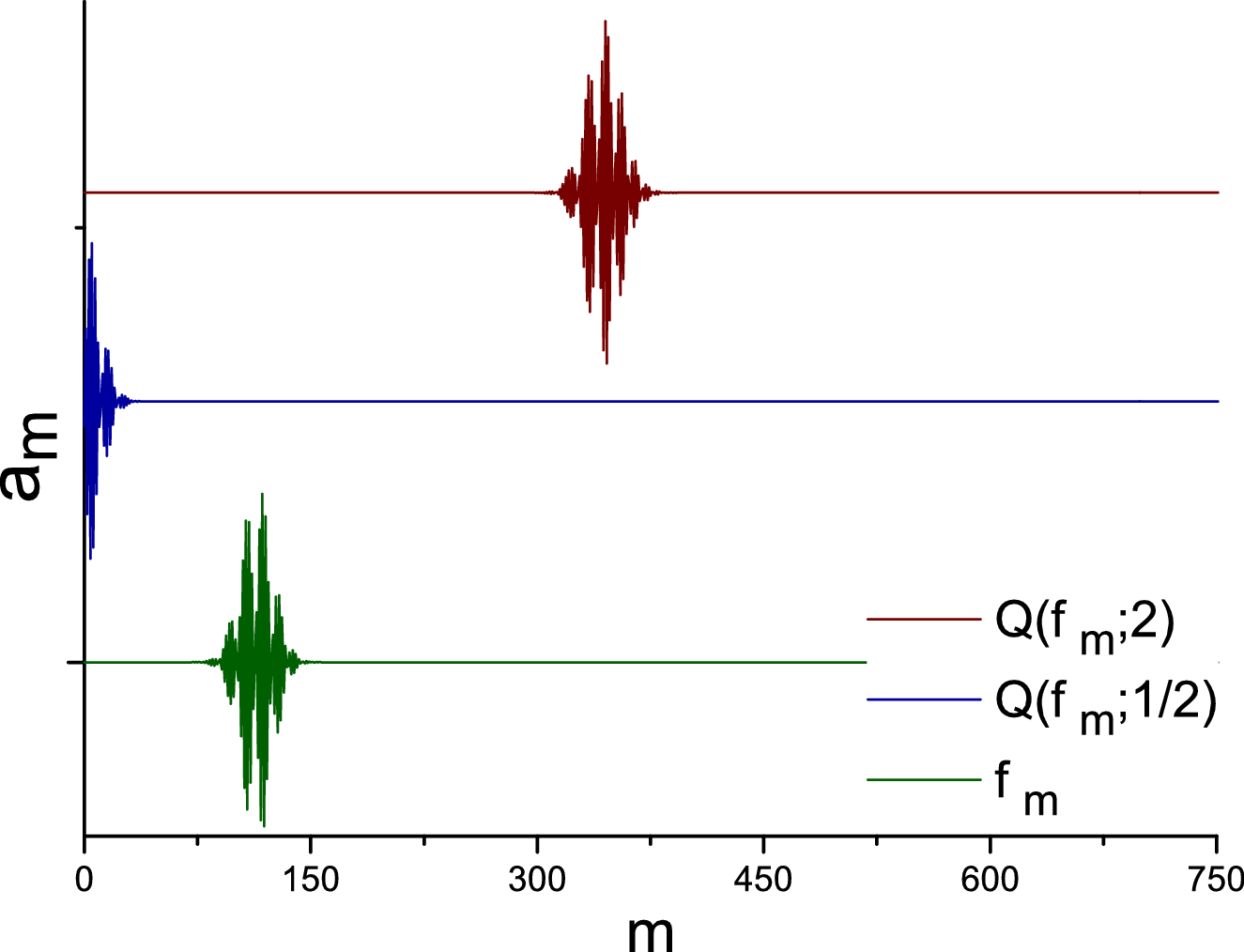}
\end{subfigure}
\caption{ a) Function (\ref{source}) and b) Laguerre spectrum for various values of parameter $\tau$ of conjugation operator $\mathbb{Q}\{\bar{f}_m;\tau\}$.}
\label{pic:reverse_test}
\end{figure}

One can see in Fig.~\ref{pic:shift_test} for formula (\ref{shift_proc}) that the expansion coefficients $\bar{g}_m=\mathbb{S}\left\{\bar{f}_m;\tau\right\}$ correspond to a function $g(t) = f(t-\tau)$, where $f(t)\equiv 0$ for $t < 0$. One can see in Fig.~\ref{pic:reverse_test} for formula (\ref{reverse_proc}) that the expansion coefficients $\bar{h}_m=\mathbb{Q}\left\{\bar{f}_m;\tau\right\}$  approximate a function $h(t)=f(\tau-t)$, where $f(t) \equiv 0$ for $t < 0$. The transform $\mathbb{S}\left\{\cdot;\tau\right\}$ will be called a shift. The transform $\mathbb{Q}\left\{\cdot;\tau\right\}$ will be called conjugation for the interval $[0,\tau]$, since the transform $\mathbb{Q}\left\{\cdot;\tau\right\}$ is an analog of complex conjugation for the coefficients of the trigonometric Fourier series. To implement the transforms (\ref{shift_proc}) and (\ref{reverse_proc}), $O(n \log n)$ operations are needed, if we use algorithms based on the fast Fourier transform \cite{Nussbaumer1982} to calculate the linear convolution (\ref{shift_proc}) and the correlation (\ref{reverse_proc}).

\subsection{Time-dependent truncation of Laguerre Series}
To remove the fictitious periodicity of the function being approximated, a procedure was developed in Section \ref{section:enetgy_truncation}. This procedure, by analyzing the Laguerre series spectra, limits the number of expansion coefficients to separate the first period of the function being approximated from all subsequent fictitious periods. Here we propose another algorithm to remove the periodicity with less computational costs without any additional increase in the approximation interval.

\begin{figure}[!htb]
\centering
\begin{subfigure}[b]{0.47\textwidth}
\includegraphics[width=\textwidth]{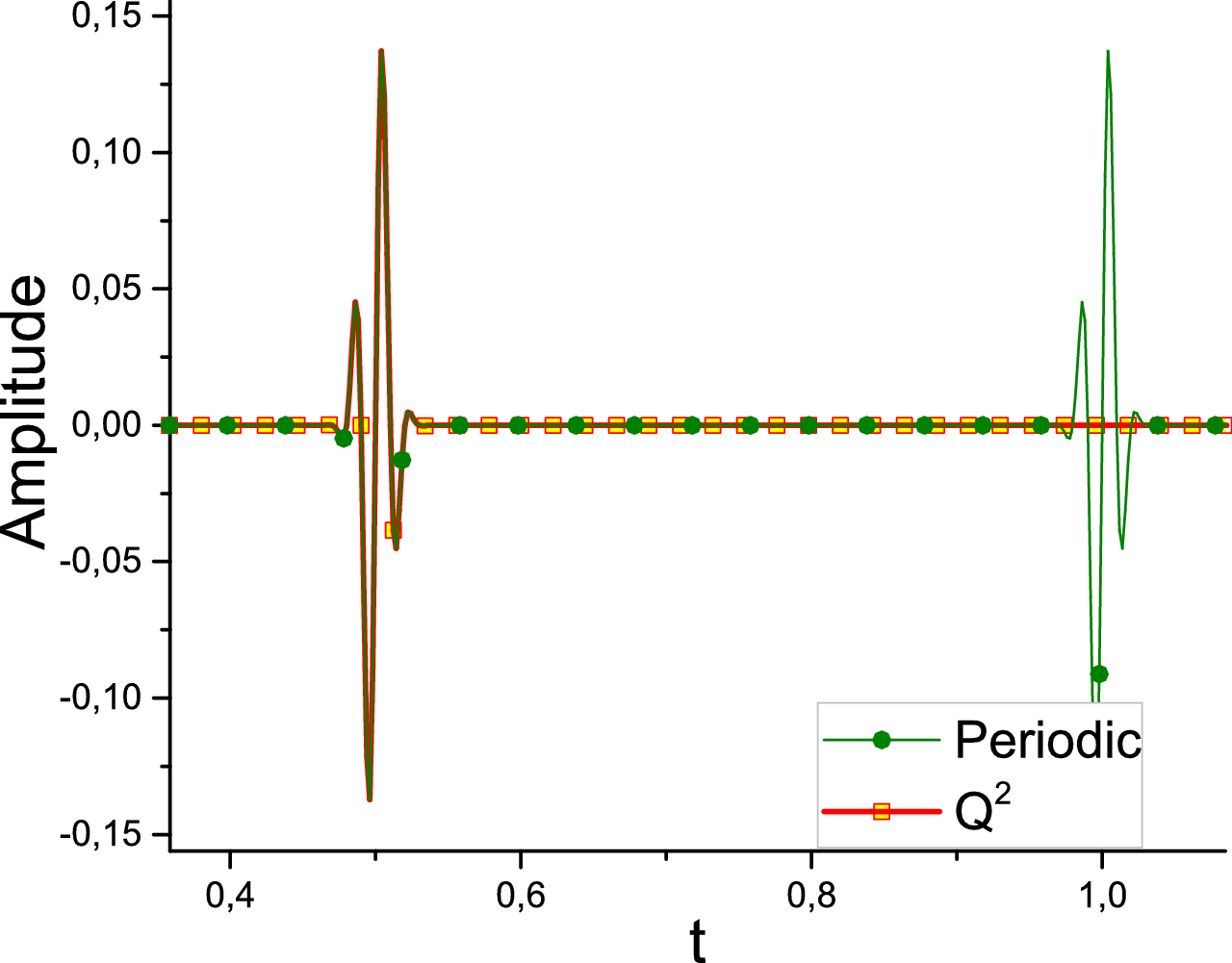}
\end{subfigure}
~
\begin{subfigure}[b]{0.47\textwidth}
\includegraphics[width=\textwidth]{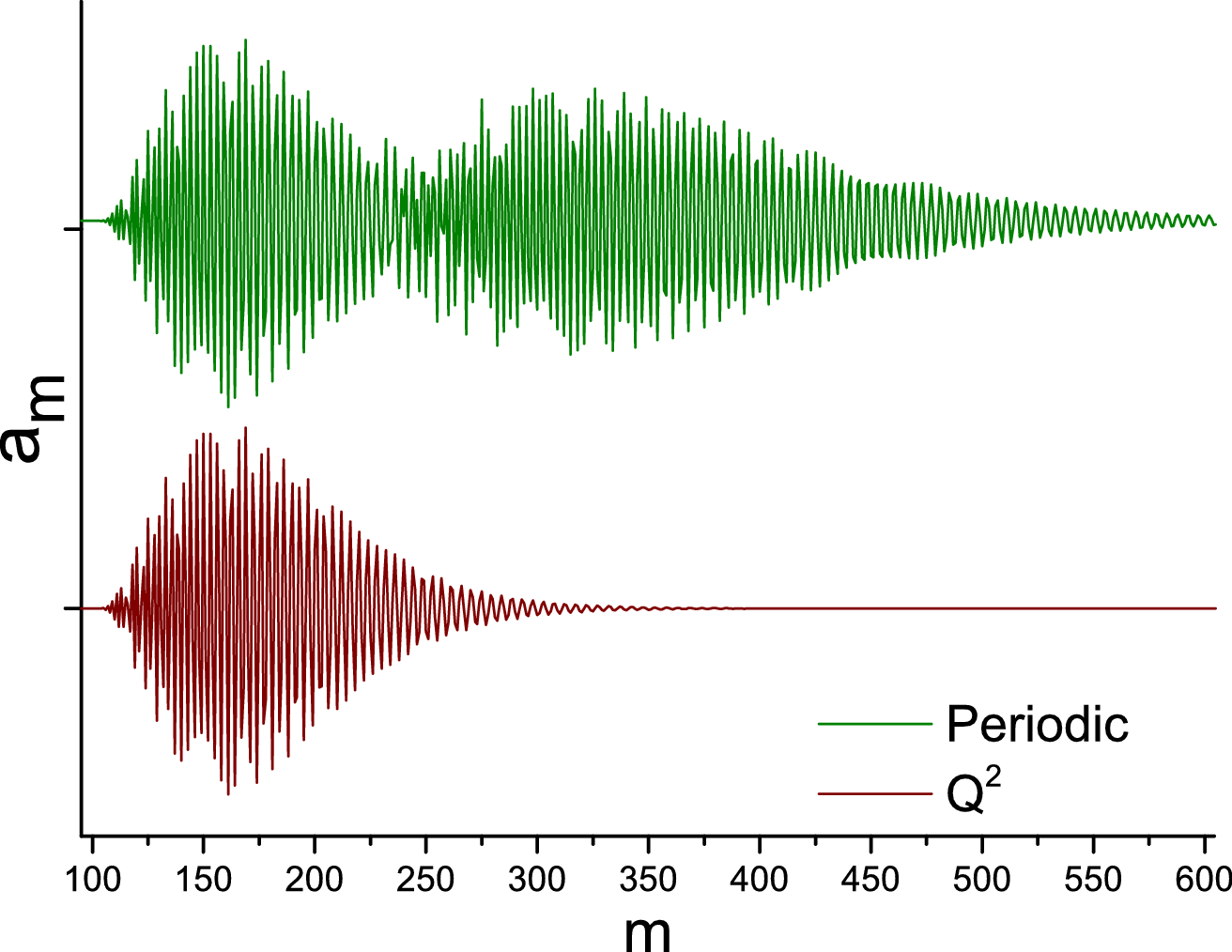}
\end{subfigure}
\caption{ a) Function (\ref{source}) before and after applying the operator $\mathbb{Q}^2\left\{\bar{f}_m;0.7\right\}$ and b) Laguerre spectrum.}
\label{pic:conjg3}
\end{figure}

Consider a procedure which, for a given parameter $\tau>0$, transforms the Laguerre series coefficients to make the series for the function $f(t)$ approximate the function $r(t)=H(-t+\tau)f(t)$, where $H(t)$ is the Heaviside function. This is equivalent to nullifying the values of the series $ \forall \  t>\tau>0$.  This can be achieved by successively applying two conjugation operations of the form $\mathbb{Q}^2\left\{\cdot;\tau\right\}\equiv\mathbb{Q}\left\{\mathbb{Q}\left\{\cdot;\tau\right\};\tau\right\}$ to the Laguerre series.
\begin{figure}[!htb]
\centering
\begin{subfigure}[b]{0.47\textwidth}
\includegraphics[width=\textwidth]{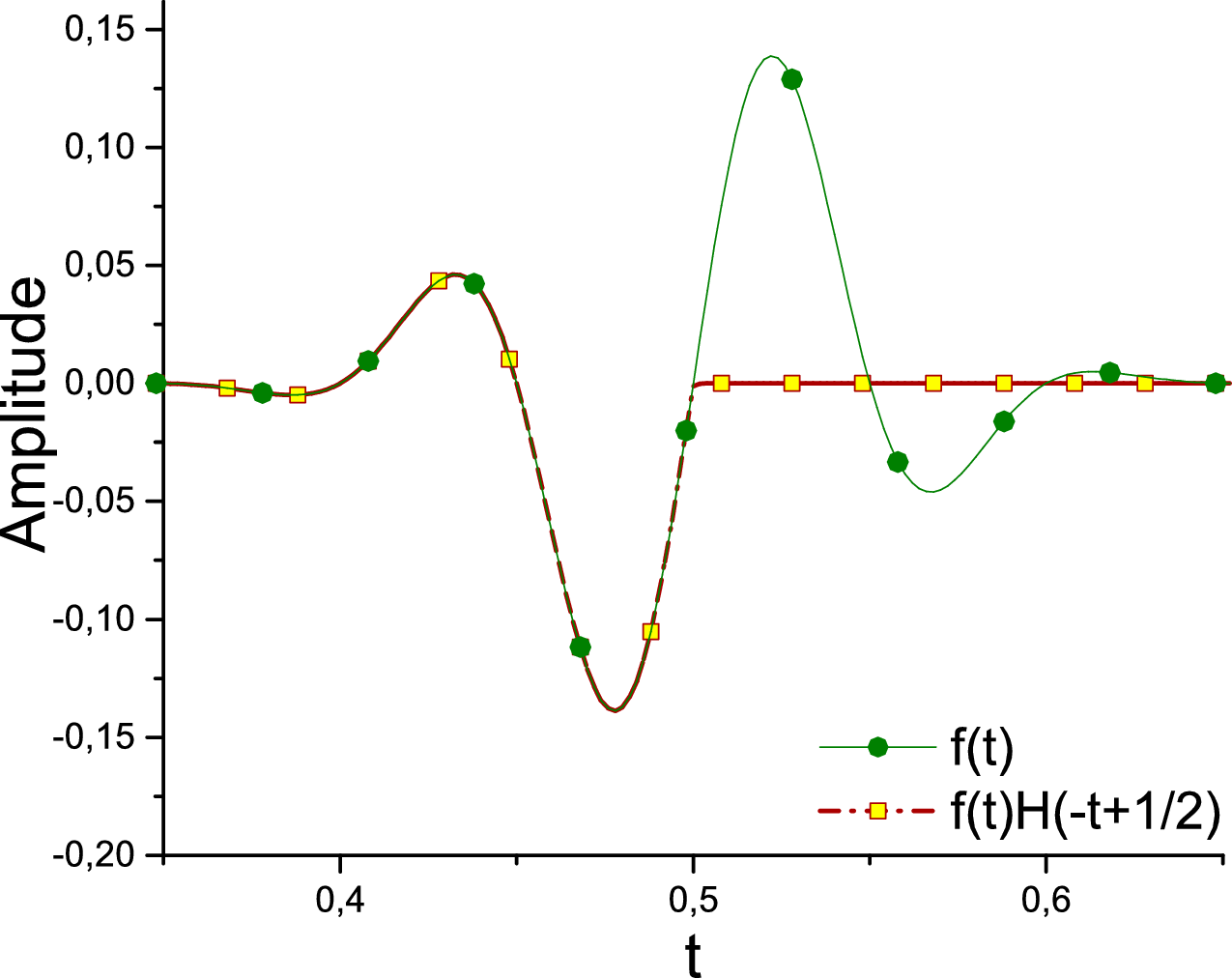}
\end{subfigure}
~
\begin{subfigure}[b]{0.47\textwidth}
\includegraphics[width=\textwidth]{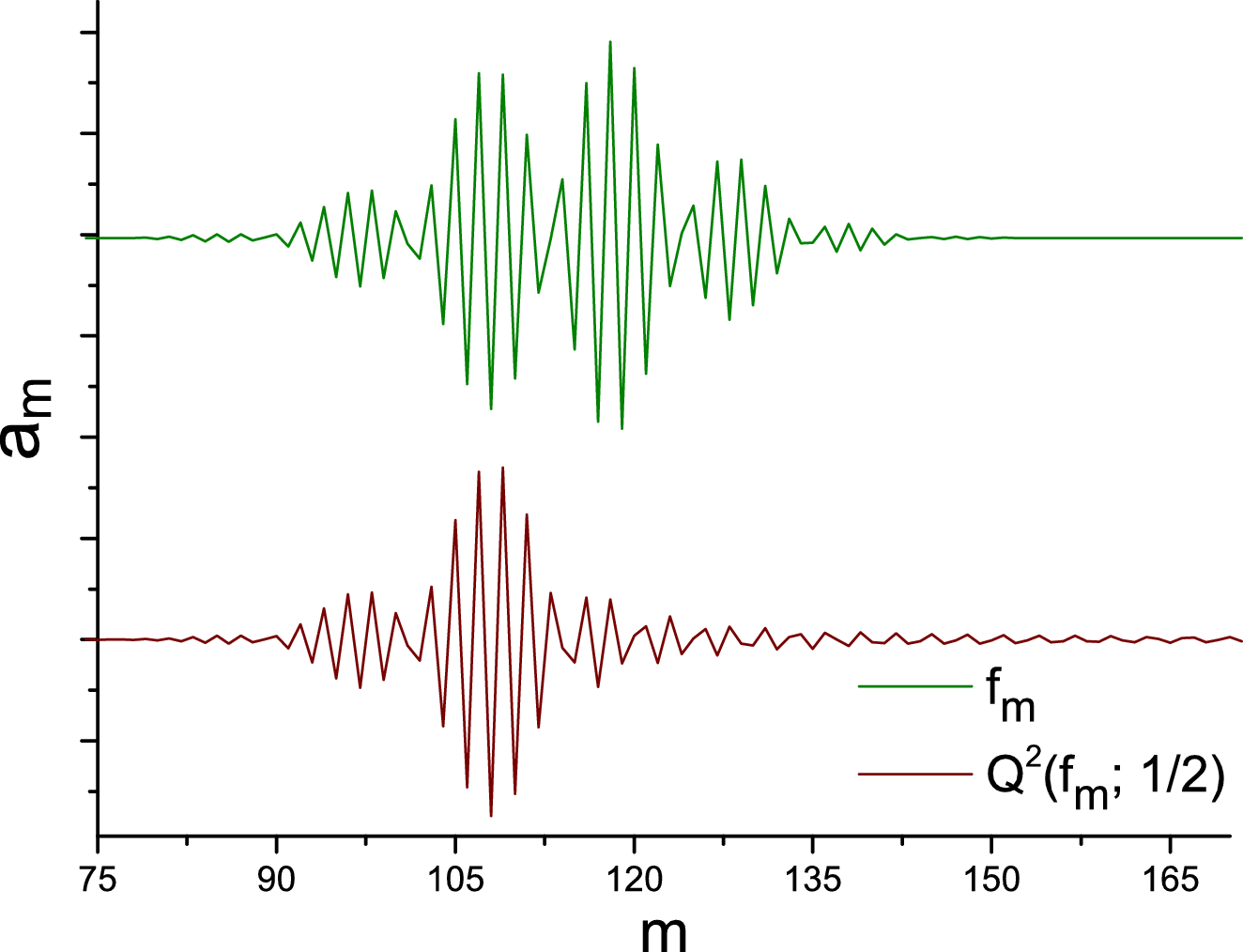}
\end{subfigure}
\caption{ Function (\ref{source}) before and after applying the operator  $\mathbb{Q}^2\left\{\bar{f}_m;1/2\right\}$ and b) Laguerre spectrum.}
\label{pic:conjg2}
\end{figure}
It follows from Fig.~\ref{pic:conjg2}a that once the operation $\mathbb{Q}^2\left\{\bar{f}_m;1/2\right\}$ is applied, the values of the series become zero,  $ \forall \  t>1/2$. The local smoothness of the function $r(t)$ in the vicinity of the point $t=1/2$ decreases, which increases the spectrum width (Fig.~\ref{pic:conjg2}b). However, if there are no additional discontinuities of the function and its derivatives at the point $t=\tau$ (Fig.~\ref{pic:conjg3}a), the spectrum width does not increase.
As a result, one can preliminarily calculate the expansion coefficients by formula  (\ref{main_formula}), and then apply the operation $\mathbb{Q}^2\left\{\cdot;T\right\}$ to remove the fictitious periodicity. The operation $\mathbb{Q}^2$ uses $O(n\log n)$ arithmetic operations, which is much less than the computational costs for formula (\ref{main_formula}). Therefore, the total costs of the approach being proposed will increase insignificantly. To avoid any additional discontinuities and decreases in the smoothness of the function being approximated and, hence, increases in the Laguerre spectrum width, the function being expanded in the series is locally multiplied by an exponentially attenuating multiplier on the right boundary of the approximation interval.

In solving practical problems of seismic prospecting, it is often necessary to perform integral transforms for a set of independent time series, called seismic traces.  In this case the procedure of removing the periodicity can be implemented in a more efficient way.  For this formula (\ref{main_formula}) is rewritten in matrix form as follows:
\begin{equation}
\label{main_matirx}
  \left(
    \begin{array}{c}
      \bar{a}_0 \\
      \bar{a}_1 \\
      ... \\
      \bar{a}_{n-1} \\
      \bar{a}_n \\
    \end{array}
  \right)=\left(    
            \begin{array}{ccccc}
              \frac{1}{\left(-{\mathrm{i}k_0+{\eta}/{2}}\right)} &\frac{1}{\left(-{\mathrm{i}k_1+{\eta}/{2}}\right)} & ... & \frac{1}{\left(-{\mathrm{i}k_{N_x}+{\eta}/{2}}\right)} \\\\
             \frac{\left(-{\mathrm{i}k_0-{\eta}/{2}}\right)}{\left(-{\mathrm{i}k_0+{\eta}/{2}}\right)^{2}} &\frac{\left(-{\mathrm{i}k_1-{\eta}/{2}}\right)}{\left(-{\mathrm{i}k_1+{\eta}/{2}}\right)^{2}} & ... & \frac{\left(-{\mathrm{i}k_{N_x}-{\eta}/{2}}\right)}{\left(-{\mathrm{i}k_{N_x}+{\eta}/{2}}\right)^{2}} \\
              ... & ... & ... & ... \\
              \frac{\left(-{\mathrm{i}k_0-{\eta}/{2}}\right)^n}{\left(-{\mathrm{i}k_0+{\eta}/{2}}\right)^{n+1}} &\frac{\left(-{\mathrm{i}k_1-{\eta}/{2}}\right)^n}{\left(-{\mathrm{i}k_1+{\eta}/{2}}\right)^{n+1}} & ... & \frac{\left(-{\mathrm{i}k_{N_x}-{\eta}/{2}}\right)^n}{\left(-{\mathrm{i}k_{N_x}+{\eta}/{2}}\right)^{n+1}} \\
            \end{array}
          \right)\left(
                         \begin{array}{c}
                           \tilde{f}_0 \\
                           \tilde{f}_1 \\
                           ... \\
                           \tilde{f}_{N_x-1} \\
                           \tilde{f}_{N_x} \\
                         \end{array}
                       \right)=M\tilde{F}.
\end{equation}
Consider the matrix $\check{M}$ whose columns are obtained from the columns of the matrix $M$ by applying the operations $\mathbb{Q}^2\left\{\cdot;T\right\}$. Instead of applying the operation $\mathbb{Q}^2\left\{\cdot;T\right\}$ to the calculated coefficients $\bar{a}_n$, one can preliminarily calculate the matrix $\check{M}$ and then calculate the expansion coefficients without the fictitious periodicity. This method is used if the number of columns of the matrix $M$ is much less than the number of functions to be approximated.

Now let us formulate algorithms to approximate a function $f(t)$, $t\in\left[0,T\right]$ by a Laguerre series, where to remove the fictitious periodicity we use the operator $\mathbb{Q}^2\left\{\cdot;T\right\}$.
\\\\
{\bf Algorithm 2} to approximate a function $f(t)$ on an interval $t\in\left[0,T\right]$ by a Laguerre series
\begin{enumerate}
  \item {\bf Preparation stage:}
  \subitem 1.1 Create a matrix $M$ of the form (\ref{main_matirx}).
  \subitem 1.2 	Calculate the modified matrix $\check{M}$ by making the transform $\mathbb{Q}^2\left\{\cdot;T\right\}$ for each column of the matrix $M$.
  \item {\bf For each of the functions $f(t)$ being approximated:}
  \subitem 2.1  	Calculate $\tilde{f}=FFT(f)$ using a fast algorithm of the discrete Fourier transform.
  \subitem 2.2  	Calculate the Laguerre series coefficients as  $\left(\bar{a}_0,\bar{a}_1,...,\bar{a}_n\right)^{T}=\hat{M}\left(\tilde{f}_0,\tilde{f}_1,...,\tilde{f}_{N_x}\right)^T$.
\end{enumerate}
If the number of functions to be approximated is smaller than the number of columns of the matrix $M$, the following algorithm, which does not calculate the matrix $\check{M}$, is more efficient:
\\\\
{\bf Algorithm 3} to approximate a function $f(t)$ on an interval $t\in\left[0,T\right]$ by a Laguerre series:
\begin{enumerate}
  \item   Calculate $\tilde{f}=FFT(f)$ using a fast algorithm of the discrete Fourier transform.
  \item   Calculate the Laguerre series coefficients as $\left(\bar{a}_0,\bar{a}_1,...,\bar{a}_n\right)^{T}=M\left(\tilde{f}_0,\tilde{f}_1,...,\tilde{f}_{N_x}\right)^T$.
  \item   Transform $\mathbb{Q}^2\left\{\bar{a}_n;T\right\}$ to exclude  the fictitious periodicity.
\end{enumerate}

\subsection{A generalization for the expansion algorithms}
\label{section:stable_laguerre_evaluation}
\begin{figure}[!htb]
\centering
\begin{subfigure}[b]{0.47\textwidth}
\includegraphics[width=\textwidth]{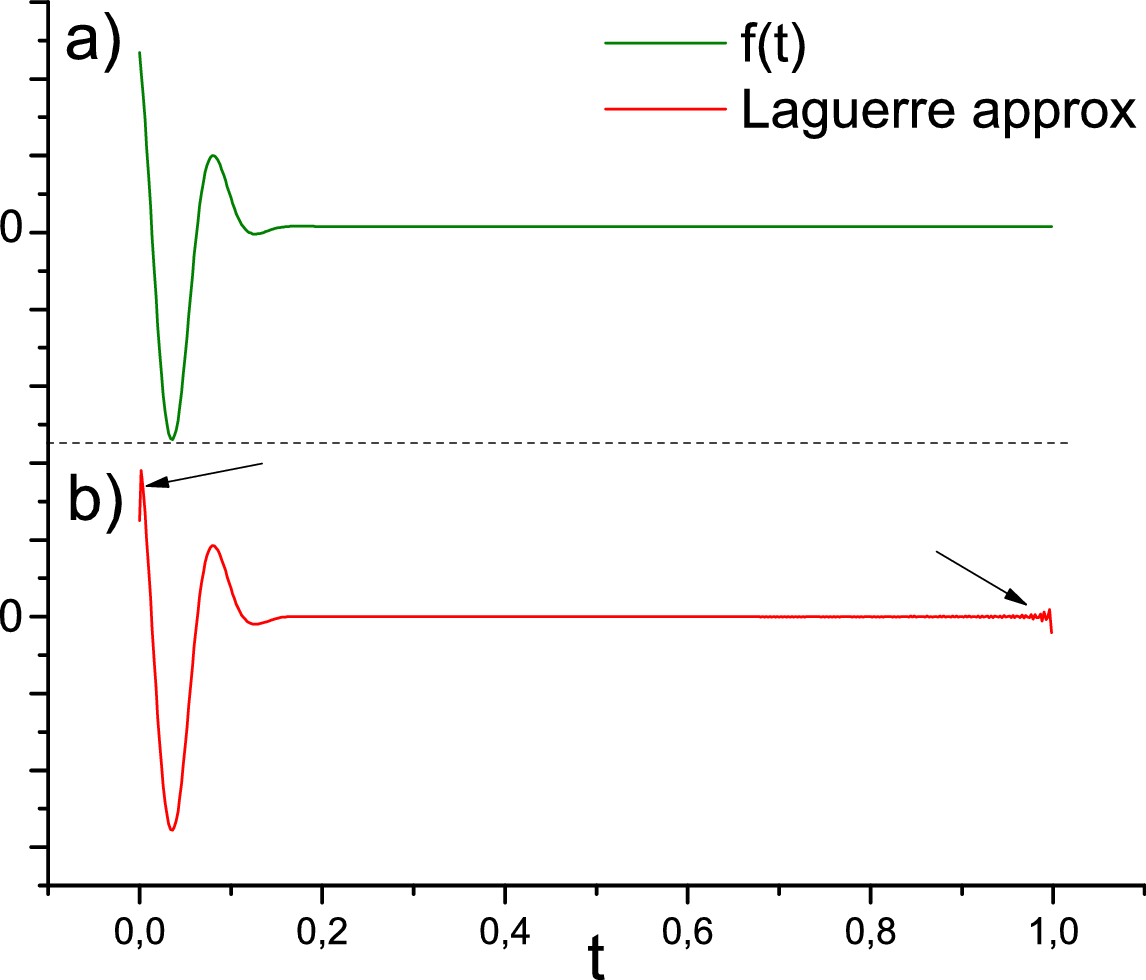}
\end{subfigure}
~
\begin{subfigure}[b]{0.47\textwidth}
\includegraphics[width=\textwidth]{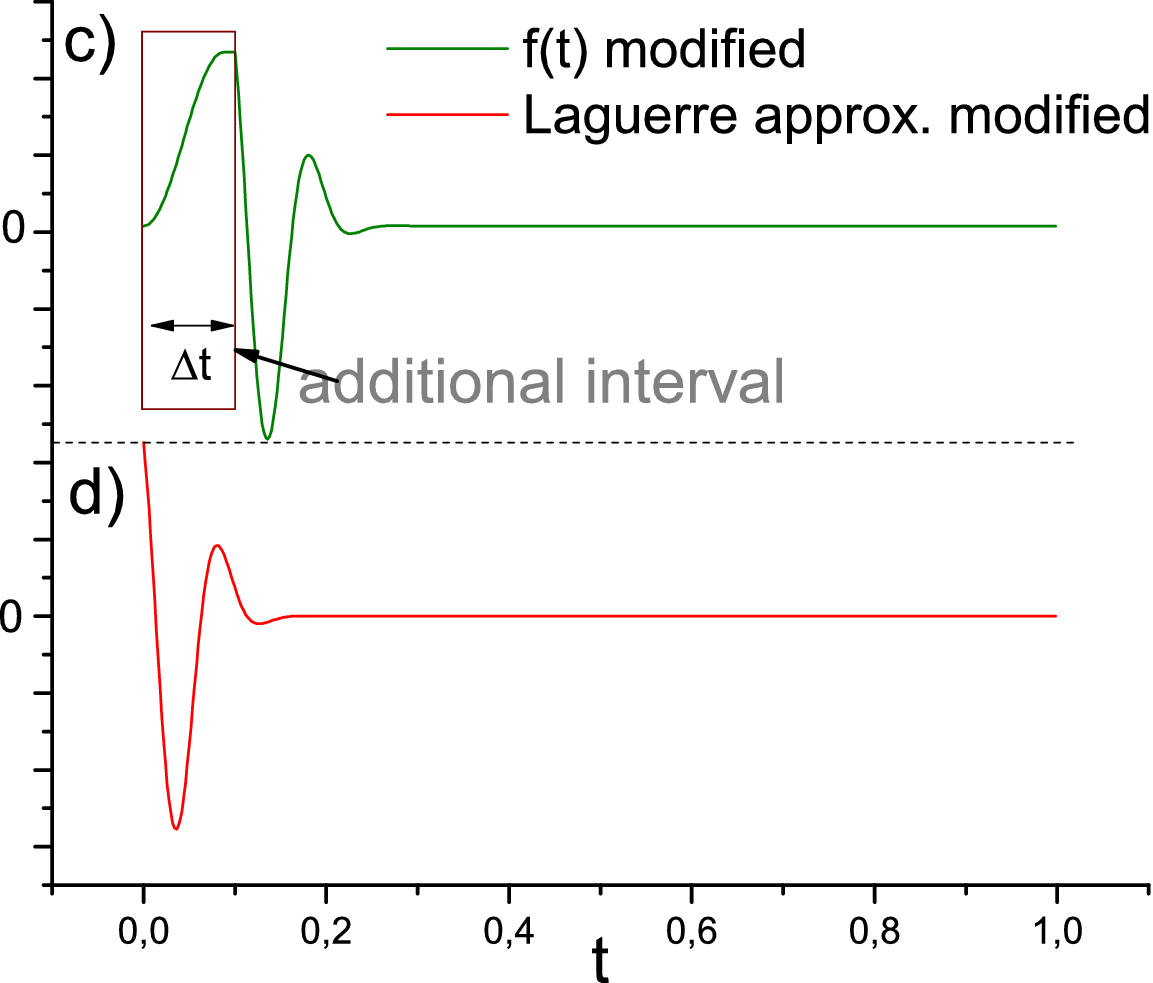}
\end{subfigure}
\caption{a) Function to be approximated, b) incorrect approximation of the initial function by a Laguerre series with artefacts shown by arrows, c) auxiliary function including an additional interval, d) correct approximation of the initial function by a Laguerre series after removing the auxiliary interval.}
\label{pic:conjg4}
\end{figure}

Algorithms 1, 2, and 3 can be used when the function to be approximated by a Laguerre series can be represented by a Fourier series as well. For the Laguerre series coefficients to decrease rapidly enough, the function being expanded must tend to zero exponentially in the vicinity of the right boundary of the approximation interval \cite{Boyd2001}. This can be achieved by locally multiplying the function by a factor of the form $\exp(-\mu t),\  \mu >0$. On the other hand, since the trigonometric interpolation is periodic, the condition $f(0) = f(T) = 0$ must be satisfied.  This imposes constraints on the form of the function being approximated. For instance, if the above-considered algorithms are applied to the function shown in Fig.~\ref{pic:conjg4}a, there will be oscillations on both boundaries of the expansion interval (see Fig.~\ref{pic:conjg4}b). The loss of accuracy can be avoided if the calculations are made by the following scheme.

First the initial function is shifted to the right by $\Delta t$. Then, on the interval $t\in[0,\Delta t]$, a smooth function taking a zero value at $t=0$ is added (see Fig.~\ref{pic:conjg4}c where a scaled quarter-period function of $\cos^2(t)$ is specified on the interval $t\in[0,\Delta t]$. If the modified function is expanded in a Laguerre series using algorithm 2 or 3, the operation $\mathbb{Q}^2\left\{\bar{a}_m;T\right\}$ is used instead of the operation $\mathbb{Q}\left\{\mathbb{Q}\left\{\bar{a}_m;T+\Delta t\right\};T\right\}$. This will make it possible to remove both the fictitious periodicity and the auxiliary interval $t\in[0,\Delta t]$. If it is planned to use algorithm 1, for which the operation $\mathbb{Q}^2\left\{\bar{a}_m;T\right\}$ is not needed, the additional interval $t\in[0,\Delta t]$ is excluded by  calculating the expansion coefficients by formula (\ref{main_formula}), but setting $\bar{a}_m=\bar{v}_m(-\Delta t)$ instead of $\bar{a}_m=\bar{v}_m(0)$.
\subsection{Stable calculation of Laguerre functions for any order and argument value}
Consider a problem of calculating Laguerre functions by performing the operations $\mathbb{S}\left\{\cdot;\tau\right\}$  and $\mathbb{Q}\left\{\cdot;\tau\right\}$. If the argument of the functions $l_m(\eta\tau)$ for (\ref{shift_proc}), (\ref{reverse_proc}) is too large, then (as noted in the introduction) there emerges an error of "overflow"{} in calculating the function $L_n(\eta\tau)$ or an error of "underflow" in calculating $\exp{(-\eta\tau/2)}$. The use of $128$-bit arithmetic does not exclude errors of these types for larger values of the argument or the order  of the Laguerre function. Therefore, we consider a more universal approach.

It follows from the relation
$
l_m(\eta t_0)=\int_{0}^{\infty}\delta(t-t_0)l_m(\eta t)dt,
$
where $\delta(t)$ is the delta function, that the coefficients $\bar{a}_k=l_m(0)=1$ of the Laguerre series (\ref{series_lag.sum}) correspond to $\delta(0)$. Then the Laguerre function can be calculated for any values of the argument using a series of shifts of the form
\begin{equation}\label{eq:shift_laguerre_evaluation}
  \{l_m(t_0)\}=\mathbb{S}\left\{...\mathbb{S}\left\{\mathbb{S}\left\{l_m(0);\tau_1\right\};\tau_2\right\}...;\tau_p\right\},\quad t_0=\sum_{i=1}^{p}\tau_i.
\end{equation}
The maximum value of the shift parameter $\tau_i$  for $64$-bit arithmetic is limited by the capacity of representing the quantity $\exp(-\eta\tau_i/4)$ for real numbers. According to the IEEE standard describing a representation of real numbers with $64$-bit precision, by choosing $\eta\tau_i\leq 4\left|\ln(2.225\times10^{-308})\right|\approx 2600$ $l_n(\eta\tau_i)$ can be calculated by (\ref{eq:little_laguerre_function})  without situations of the "underflow"{} or "overflow"{} type.
\begin{figure}[!htb]
\centering
\begin{subfigure}[b]{0.47\textwidth}
\includegraphics[width=\textwidth]{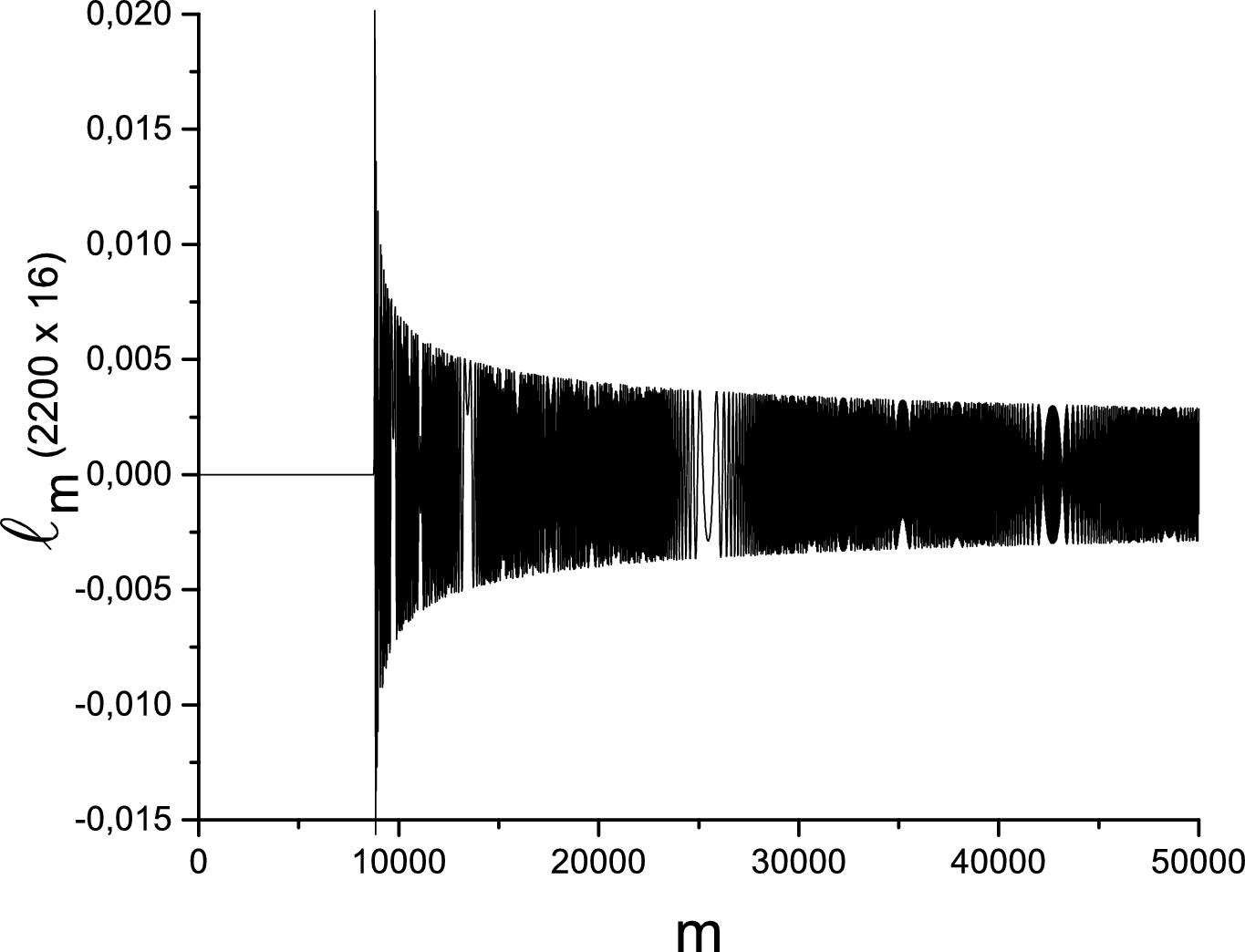}
\end{subfigure}
~
\begin{subfigure}[b]{0.47\textwidth}
\includegraphics[width=\textwidth]{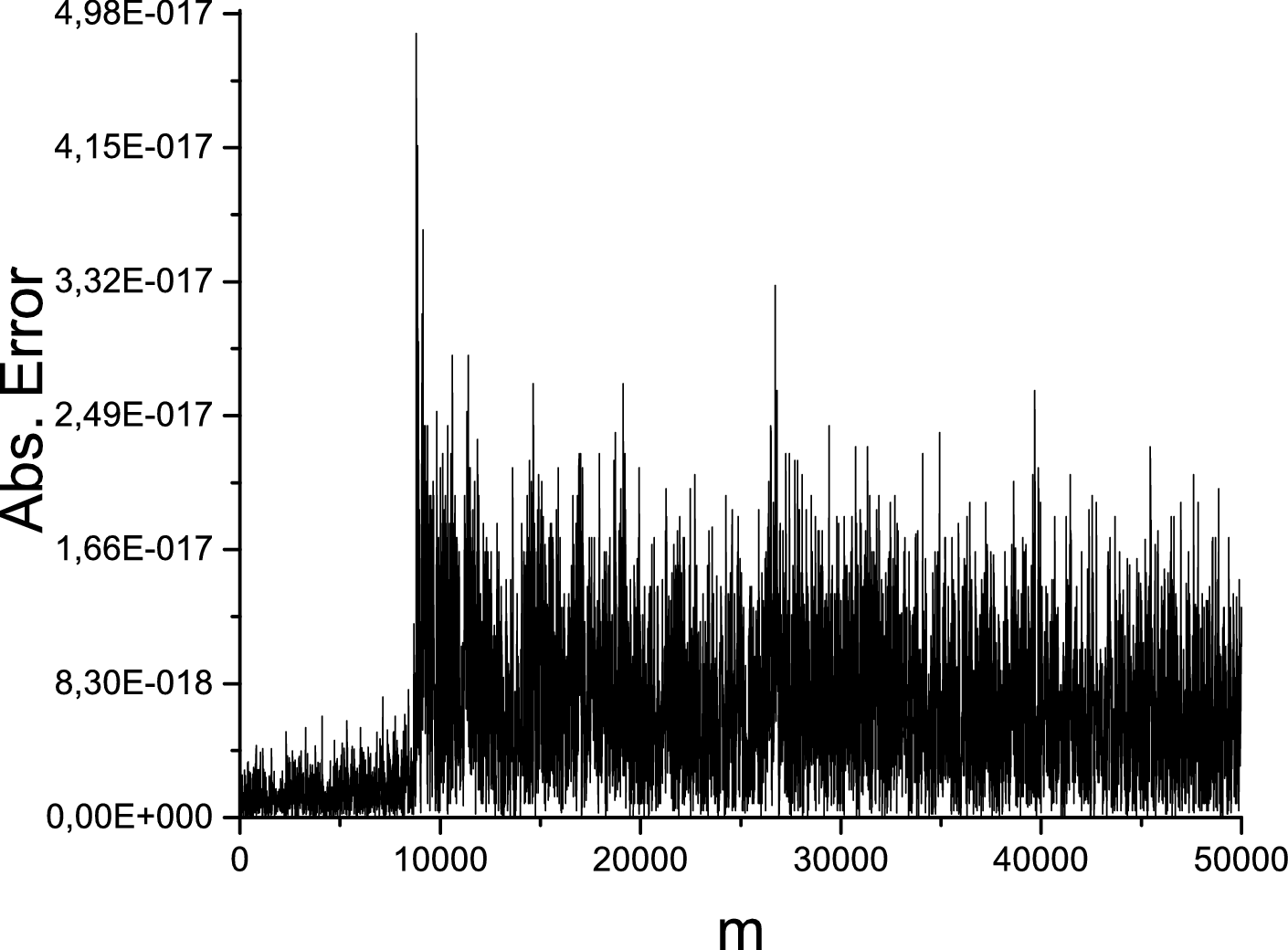}
\end{subfigure}
\caption{a) Function $l_m(\eta t)$ versus m for a constant value of argument $\eta t=2000\times 16$, b) difference of Laguerre function values calculated by formula (\ref{eq:reccurence_shift}) in $64$-bit arithmetic and formula (\ref{eq:little_laguerre_function}) in $128$-bit arithmetic.}
\label{pic:conjg24}
\end{figure}

To decrease the total number of shifts and, hence, the computational costs, it makes sense to perform the shifts recurrently:
\begin{equation}
\{l_m(2^p\eta\tau)\}=\mathbb{S}\left\{\mathbb{S}\left\{\mathbb{S}\left\{\mathbb{S}\left\{l_m(\eta\tau);\tau\right\};2\tau\right\},4\tau\right\}...;2^{p-1}\tau\right\}.
\label{eq:reccurence_shift}
\end{equation}
In comparison to formula(\ref{eq:shift_laguerre_evaluation}), the number of calculations can be reduced owing to the fact that the Laguerre function values obtained at the previous step are used in formula (\ref{shift_proc}) to make the shift at the current step of implementing formula (\ref{eq:reccurence_shift}). Thus, the shift value at each step is doubled, which decreases the total number of shifts with each of them requiring $O(n\log n )$ operations. Note that for the first shift $l_m(\eta\tau)$ must always be calculated by formulas (\ref{eq:little_laguerre_function}) and (\ref{recurrence_laguerre}).

Fig.~\ref{pic:conjg24}a shows the result of calculation of the functions $l_m(2200\times16)$ by formula (\ref{eq:reccurence_shift}). The absolute difference of the values for formula (\ref{eq:reccurence_shift}) in $64$-bit arithmetic and formulas  (\ref{eq:little_laguerre_function}), (\ref{recurrence_laguerre}) in $128$- bit arithmetic is shown in Fig.~\ref{pic:conjg24}b. It is evident from this figure that both approaches give practically the same results. Thus, the above algorithm does not use high-precision arithmetic in performing stable calculations of Laguerre functions of any order for any values of the argument. Moreover, some test calculations have shown that, in comparison to $128$-bit arithmetic, the above calculation method needs several times less calculation time if, in particular, $32$-bit arithmetic is used to organize the shift procedure.
\subsection{Optimization for a large  interval approximation}
%
It is well-known that, owing to the high performance and stability of the algorithm of fast Fourier transform, it has been widely used in many branches of computational mathematics, whereas  no algorithm for the Laguerre  transform having comparable efficiency has been developed so far. Although general methods of fast polynomial transforms were proposed long ago \cite{GOHBERG1994411}, they are of theoretical rather than practical importance. This is because they use numerically unstable efficient procedures of multiplying matrices $V$ and $V^T$  by a vector, where $V$ is an ill-conditioned Vandermonde matrix \cite{Pan2016,Gautschi2011}. For instance, fast multiplication by the matrix $V$ can be performed by using an algorithm \cite{Borodin1974} whose computational complexity is of the order of $O(n\log^2n)$ operations.  Unfortunately, this method is unstable, since one of its stages includes a recursive use of the operation of polynomial division. Multiplication of the matrix $V^T$ by a vector can be reduced to solving systems of linear algebraic equations with a Vandermonde matrix with an operation count of the order of $O(n\log^2n)$ \cite{GOHBERG1994411,Gohberg1994,Pan1993}. This approach also cannot be recommended for practical use due to its numerical instability.

The condition number for Laguerre functions is greater than that for the other classical orthogonal polynomials \cite{Gautschi1983}. Therefore, the problem of stability of fast algorithms for the Laguerre transform is probably one of the most difficult ones. By now, fast transforms have been developed for Chebyshev, Legendre, and Hermite polynomials \cite{Alpert1991,Hale2016,Leibon2008}. In these cases the arithmetic complexity of the algorithms is of the order of $O(n\log n)$ or $O(n \log^2 n)$ operations. Fast algorithms of changing from one orthogonal polynomial basis specified by a three-term recurrence relation to another one have also been developed \cite{Bostan2010}. In paper \cite{ONeil2010}, an algorithm for fast polynomial transforms based on an approximate factorization of the matrices $V$ or $V^T$ was proposed. In some cases the authors managed to decrease the computational costs to a level of $O(n \log n)$  arithmetic operations. However, the computational complexity may vary widely for various orthogonal polynomials and expansion interval lengths. Also, the algorithm becomes efficient in comparison to the direct method of multiplying a matrix by a vector, for $n\ge n_0$, where $n_0$ is of the order of several thousand.

\begin{figure}[!htb]
\centering
\includegraphics[width=\textwidth]{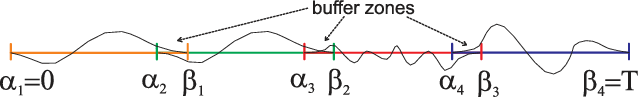}
\caption{Decompositions of the initial approximation interval into four overlapping subintervals.}
\label{pic:fast1}
\end{figure}

\begin{figure}[!htb]
\centering
\includegraphics[width=\textwidth]{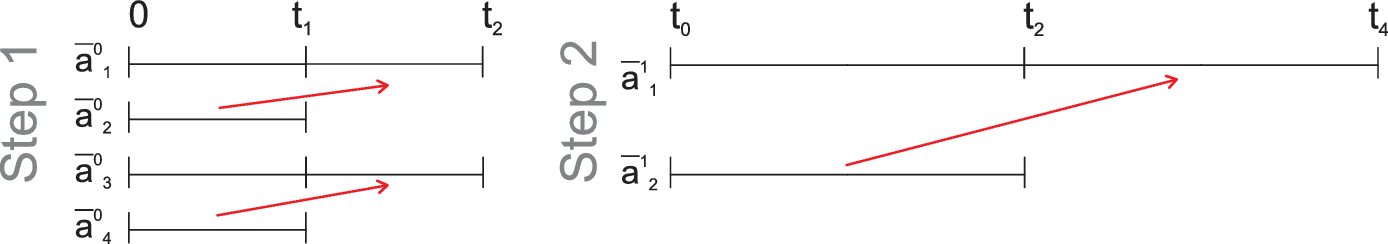}
\caption{Approximation construction scheme for function $f(t)$ with  precalculated approximations for local functions $f_i(t)$.}
\label{pic:fast2}
\end{figure}

To expand a function into a Laguerre series by formula  (\ref{main_formula}), about $O(n N_x)$ arithmetic operations are needed, where n is the number of expansion terms of the Laguerre series and $N_x$  is the number of harmonics of the auxiliary Fourier series. Approximation of the function for longer intervals calls for specifying larger values of $n$ and $N_x$, which makes the Laguerre transform inefficient. To decrease the calculation time when performing the Laguerre transform, we consider an algorithm of the "divide and conquer"{} type \cite{Smith1985}. The general idea of this approach is that at the first stage the initial problem is divided into independent subproblems with much less computational costs needed for their solution. At the second stage the solution to the initial problem is assembled from the solutions to the subproblems. This approach was successfully used, for instance, in papers in which a parallel dichotomy algorithm was proposed to solve systems of linear algebraic equations with three-diagonal \cite{terekhov:Dichotomy}, block-diagonal \cite{Terekhov:2013}, and Toeplitz matrices \cite{Terekhov2016}.
\\\\
{\bf Algorithm 4.} to approximate a function  $f(t)$ on an interval $t\in\left[0,T\right]$  by a Laguerre series:
\begin{enumerate}
  \item Decompose the approximation interval $t \in[0,T]$ into $p=2^s$ overlapping subintervals of lengths $\Delta t_i=\beta_i-\alpha_{i}$  (Fig.~\ref{pic:fast1}). In this case the function must smoothly tend to zero on the subinterval boundaries in the buffer zones so that the sum of the two local functions remains equal to the value of the function being approximated.

  \item The local function $f_i(t)$ specified on the subinterval with number $i$ is expanded in a Laguerre series on the auxiliary interval $[0,\Delta t_i]$ by algorithm 1, 2, or 3.
  \item  Shift the local functions by changing from the interval $[0,\Delta t_i]$ to the subinterval $[\alpha_i,\beta_i]$. This is done by a series of shifts of the function $f_i(t)$ using the scheme presented in Fig.~\ref{pic:fast2}, which gives an example of four subintervals. The process of assembly consists of $\log_2p$ steps, where $p$ is the number of subintervals. Hence, two steps will be needed for the example being considered. At the first step the sequences of the Laguerre series coefficients for the local functions $f_2(t)$ and $f_4(t)$ are supplemented  by zeroes to double the number of expansion coefficients.

      Then the thus expanded series are shifted using the procedures $\mathbb{S}\left\{\bar{a}_{n/2};\alpha_2\right\}$  and      $\mathbb{S}\left\{\bar{a}_{n/2};\alpha_4-\alpha_3\right\}$.  After this the corresponding coefficients of the first and second series and of the third and fourth series are added pairwise. This results in two intervals of larger lengths. At the second step this process is used for the new second series, and after it is shifted by $\mathbb{S}\left\{\bar{a}_n;\alpha_3\right\}$ the expansion coefficients of the first and second series are added. Thus, all local functions will be shifted to their initial positions with respect to the variable $t$, and the thus obtained series will approximate the initial function $f(t)$ with some accuracy.
\\{\bf Remark.} To execute one shift $\mathbb{S}\{\bar{a}_n;\tau\}$  using the fast Fourier transform, $O(n\log n)$ arithmetic operations are needed. One can see in Fig.~\ref{pic:shift_test} that a shift of the function to the right increases the number of coefficients of the Laguerre series needed to approximate the shifted function with the previous accuracy. For the calculation scheme in Fig.~\ref{pic:fast2} every shift will double the minimum number of the Laguerre series terms. Therefore, before making a shift the sequence of coefficients of the Laguerre series must be added by zeros (zero padding). After making the shift the zero values of the added expansion coefficients will become nonzero ones.
 \end{enumerate}

For larger values of $n$ and $N_x$ the computational complexity of algorithm 4 will be of the order $O(nN_x/p +n\log_2 n\log_2 p)$ vs. $O(nN_x)$, where $p$ is the number of subintervals. The first term  is the costs to approximate the local functions $f_i(t),\ t \in [0,\Delta t_i]$, and the second one is the costs to perform a series of shift operations to transform the local expansion coefficients to the expansion coefficients for the initial function $f(t)$. However, this algorithm has the following shortcoming: the number of Laguerre series coefficients to approximate the local functions $f_i(t)$ on the subintervals $[0,\Delta t_i]$ depends not only on the lengths  of the subintervals, but also on the smoothness of the functions $f_i(t)$. Taking into account that in solving practical problems the function to be approximated may have low smoothness, the convergence of the series may be not high. This results in the fact that at the same accuracy the total number of expansion coefficients for the local problems for algorithm 4 will be greater than the number of expansion coefficients when using algorithm 1, 2, or 3. Thus, the division will require additional computational costs, which can be estimated in computational experiments.
\section{Computational experiments} To estimate the accuracy of the approximation and the efficiency of the methods being proposed, let us perform a series of computational experiments to approximate functions of various smoothness on intervals of various lengths. The numerical procedures to calculate the Laguerre coefficients will be performed with single and double precision. Algorithms 2 and 3 give the same results in calculating the Laguerre series coefficients and, therefore, no separate testing of algorithm 3 will be considered.

\subsection{Inversion of Laguerre transform }
Consider the problem of calculating the inverse Laguerre transform (\ref{series_lag.sum}). In contrast to the direct transform, in the summation of the series there only remains the problem of calculating the Laguerre functions of high orders for larger argument values. This problem can be solved in several ways. If the calculations are made with $128$-bit real precision by formula (\ref{eq:little_laguerre_function}), the Laguerre functions of high orders can be calculated for rather large values of the argument without errors of "overflow"{} and "underflow"{}. Another method is to use asymptotic expansions \cite{Temme1990,Gil2017} to calculate the Laguerre polynomials $L_n(\eta t)$, whence, multiplying by $\exp{(-\eta t/2)}$, we obtain the  Laguerre functions.

If a function approximated by a Laguerre series can be represented by a Fourier series, one can change from the Laguerre coefficients to Fourier coefficients using the following formula:
\begin{equation}
\label{matrix_inv}
(\tilde{f}_0,\tilde{f}_1,...,\tilde{f}_{N_x})^T=\frac{1}{T}M^{*}(\bar{a}_0,\bar{a}_1,...,\bar{a}_n)^T,
\end{equation}
where $\tilde{M}$ is a modified matrix of the form (\ref{main_matirx}). In this case a major problem is in the emergence of discontinuities of the function on the boundaries of the approximation interval, $t\in[0,T]$.

Finally, we can use the stable method of calculating the Laguerre functions by formulas  (\ref{eq:shift_laguerre_evaluation}) or (\ref{eq:reccurence_shift}) considered in Section \ref{section:stable_laguerre_evaluation}. As noted above, although the fast Fourier transform is needed to calculate the linear convolution, this method of organizing the calculations requires less calculation time than when using $128$-bit arithmetic and formula (\ref{eq:little_laguerre_function}).
\subsection{Test 1. Expansion of a smooth function}
As a first test, consider, on the interval $t \in [0,1]$, an approximation of a function $f(t)$ of the form (\ref{source}) with parameters $f_0=30,\ g=4,\ t_0=0.5$. The discretization step of the function $ht=0.002$. The approximation error is estimated by the formula
\begin{equation}
\displaystyle \epsilon=\sqrt{\frac{\sum_{i=1}^{s}\left(f(t_i)-\sum_{j=1}^{n}\bar{a}_j l_j(\eta t_i)\right)^2}{\sum_{i=1}^{s}f^2(t_i)}},
\end{equation}
where $f(t_i)$  is the function to be expanded in a Laguerre series, which is specified on a set of values $ t_i \in [0,T], i=1,2,...s$, $t_1=0,\ t_s=T $.

\begin{figure}[!htb]
\centering
\begin{subfigure}[b]{0.47\textwidth}
\includegraphics[width=\textwidth]{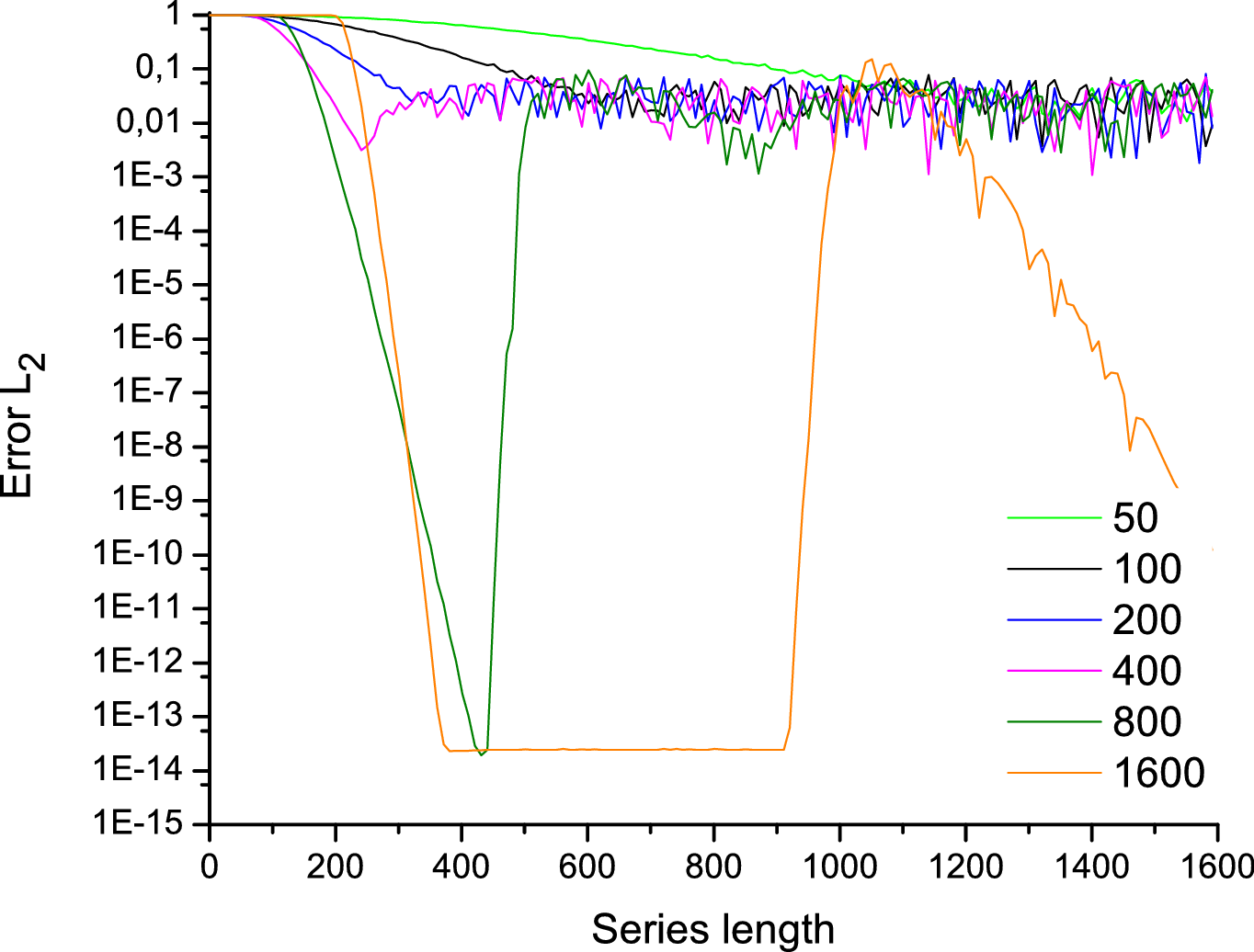}
\end{subfigure}
~
\begin{subfigure}[b]{0.47\textwidth}
\includegraphics[width=\textwidth]{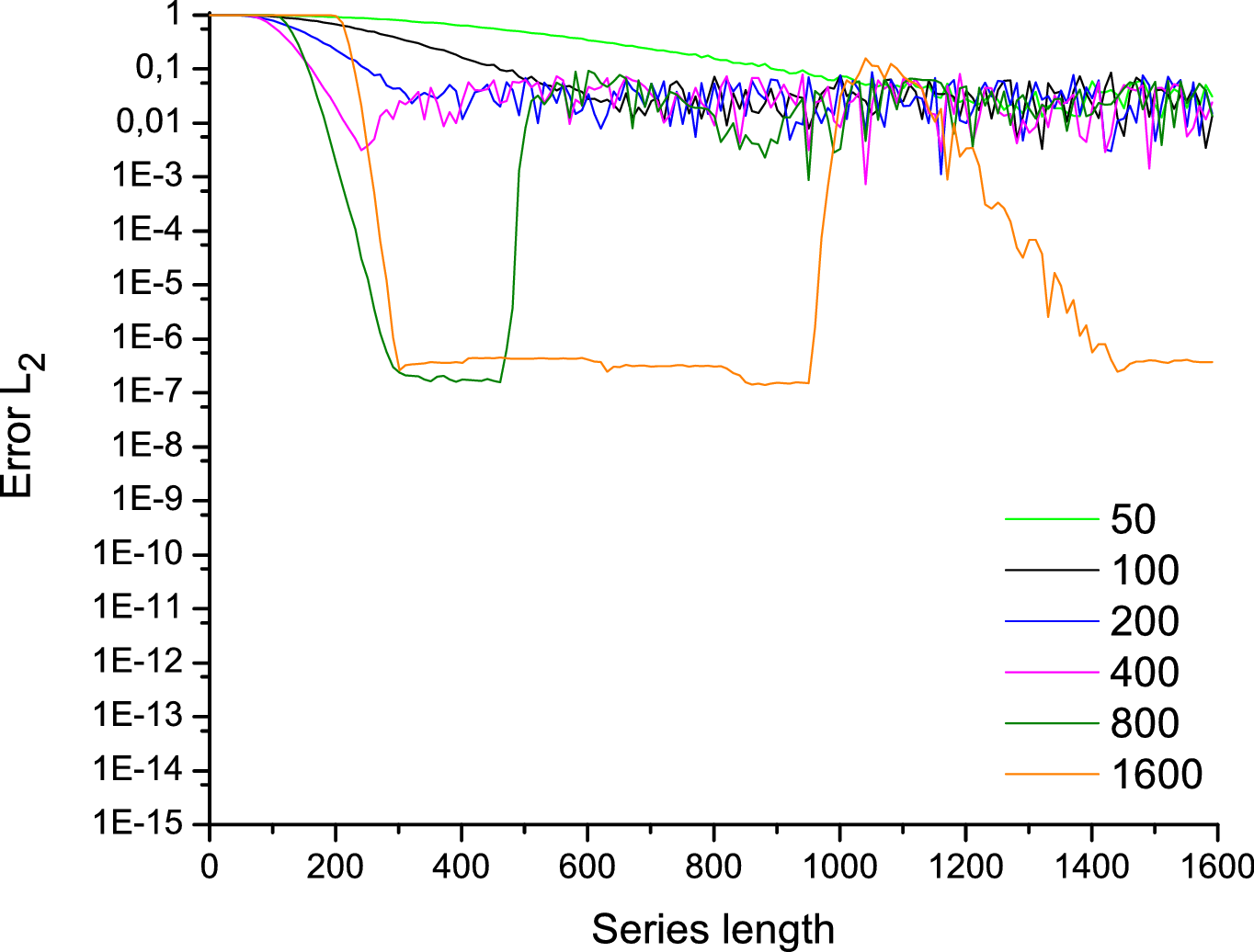}
\end{subfigure}

\begin{subfigure}[b]{0.47\textwidth}
\includegraphics[width=\textwidth]{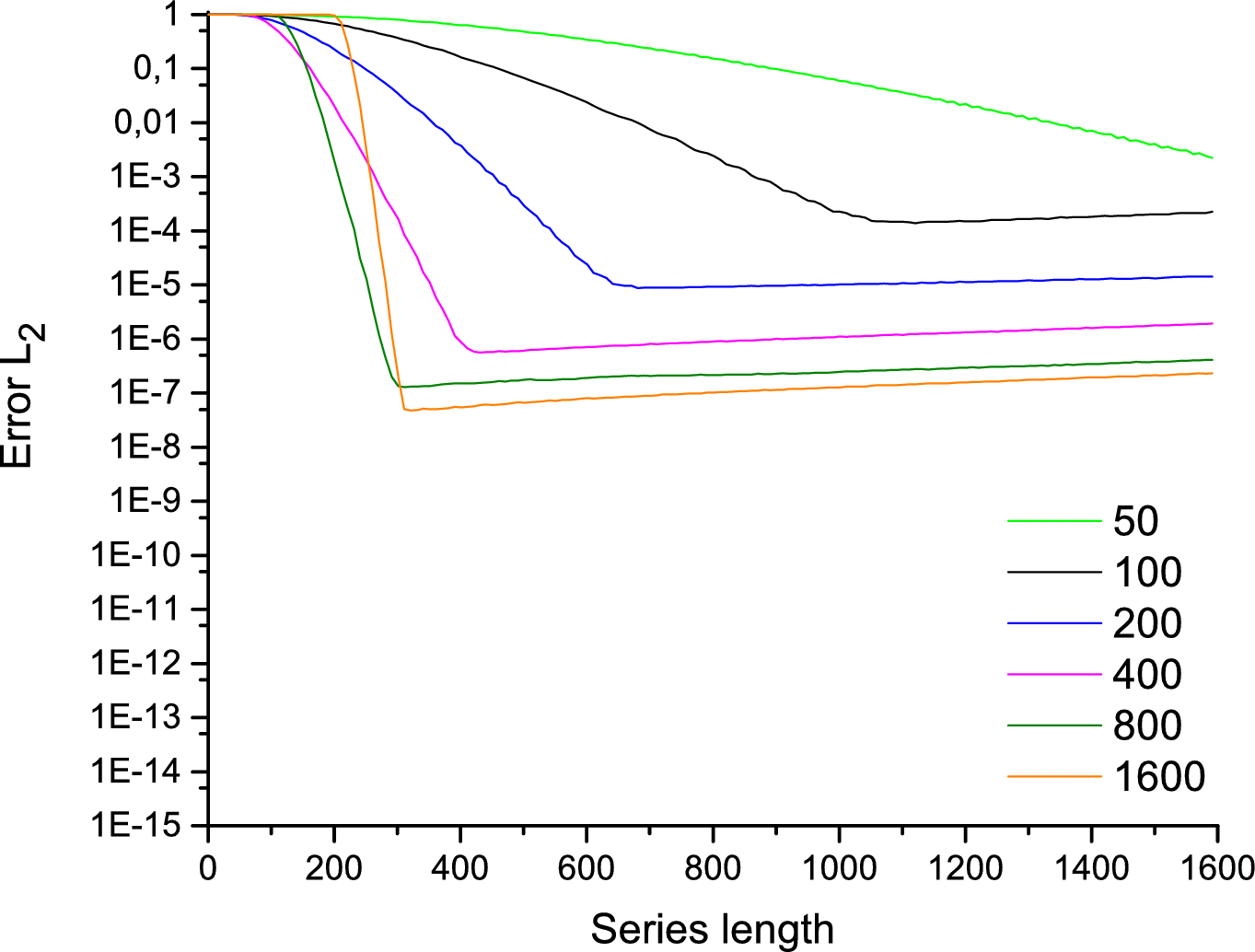}
\end{subfigure}
~
\begin{subfigure}[b]{0.47\textwidth}
\includegraphics[width=\textwidth]{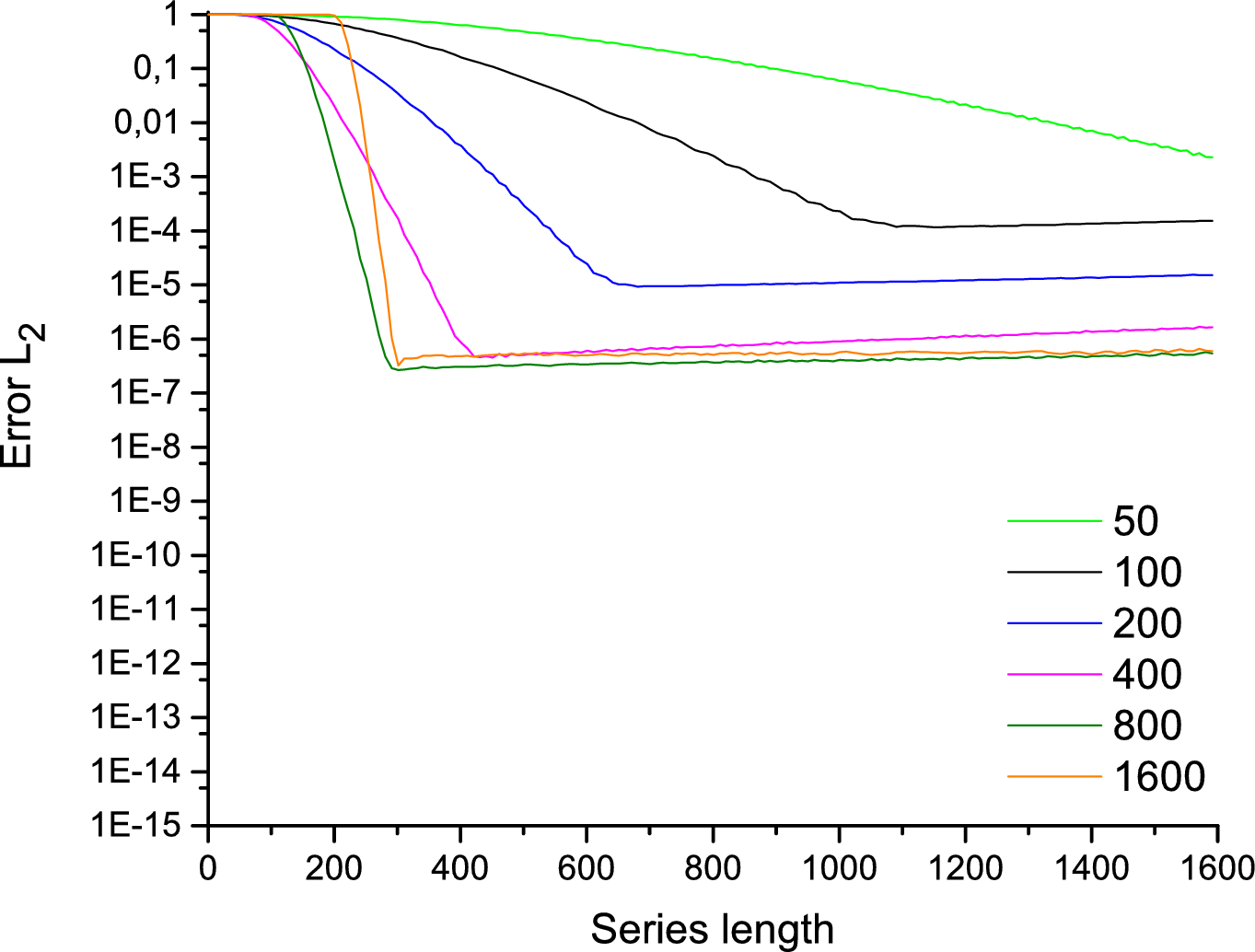}
\end{subfigure}

\caption{Approximation error for function (\ref{source})  versus the number of terms of Laguerre series for various values of the transform parameter, $\eta=50,100,...,1600$, a) algorithm 1 with 64-bit precision, b) algorithm 1 with 32-bit precision, c) algorithm 2 with 64-bit precision, d) algorithm 2 with 32-bit precision.}
\label{pic:test1}
\end{figure}

\begin{figure}[!htb]
\centering
\begin{subfigure}[b]{0.47\textwidth}
\includegraphics[width=\textwidth]{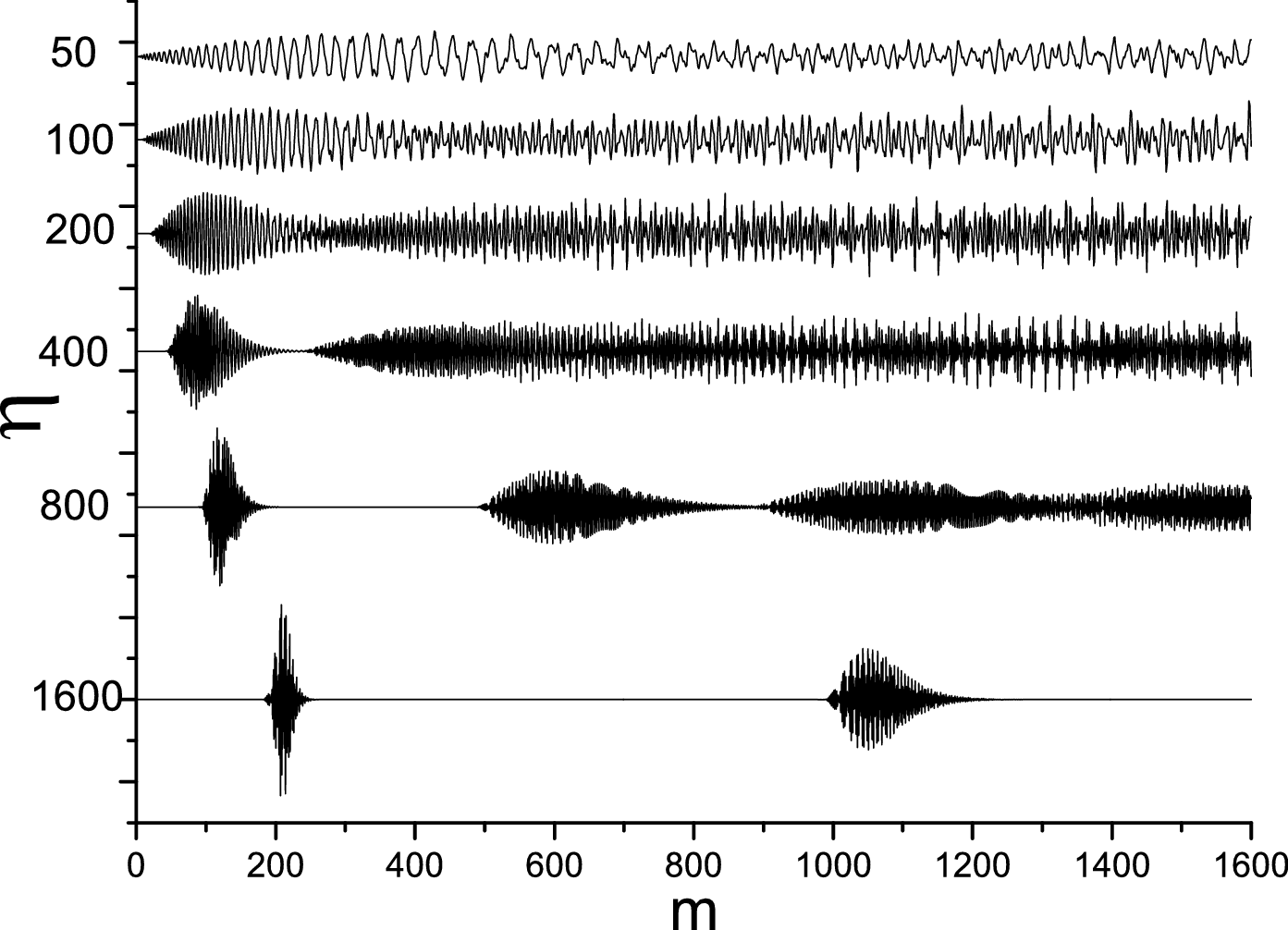}
\end{subfigure}
~
\begin{subfigure}[b]{0.47\textwidth}
\includegraphics[width=\textwidth]{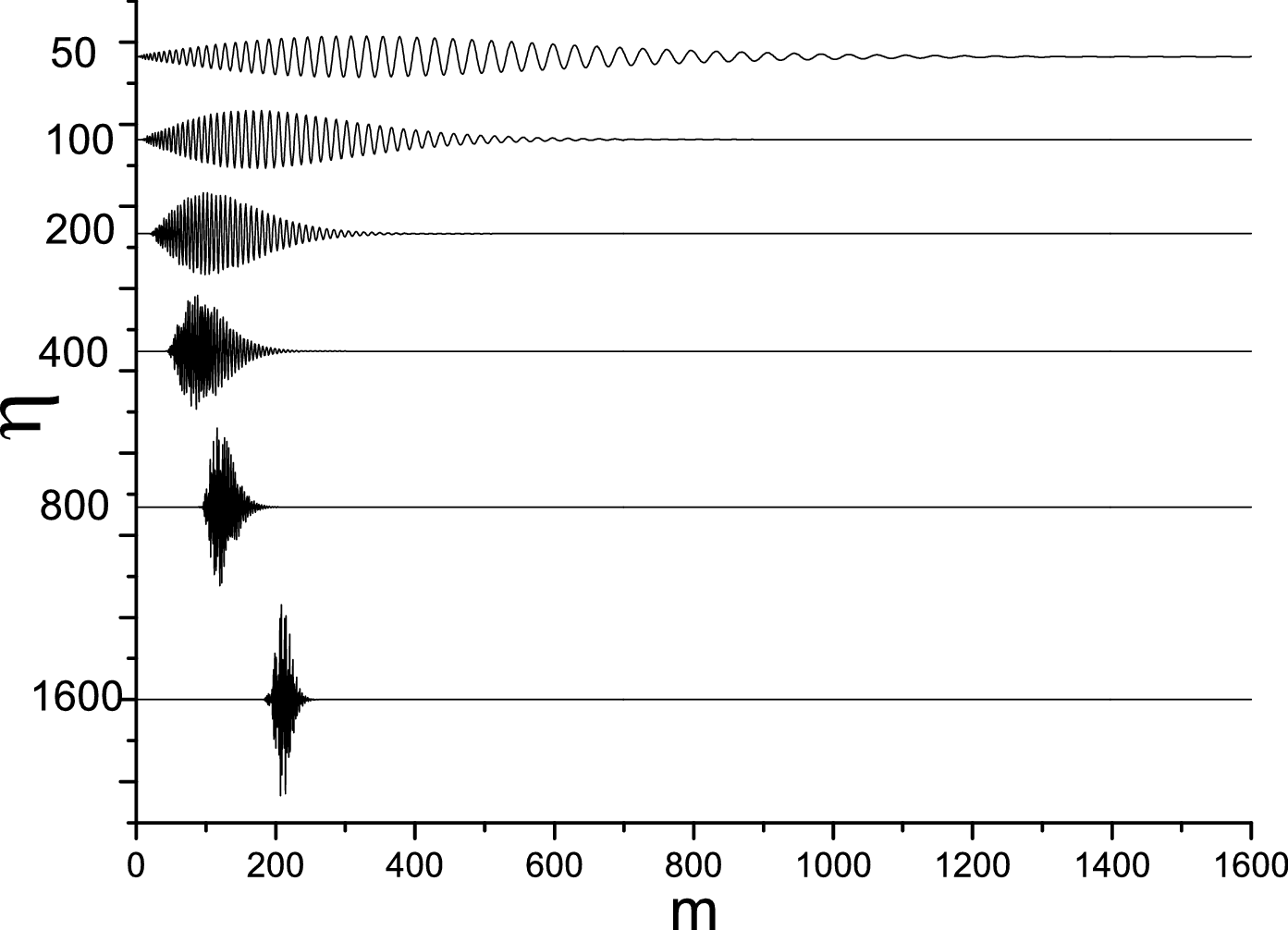}
\end{subfigure}
\caption{Laguerre spectrum for various values of the transform parameter, $\eta=50,100,...,1600$ for function (\ref{source}) for a) algorithm 1 and b) algorithm 2.}
\label{pic:test11}
\end{figure}

Fig.~\ref{pic:test1} shows the error versus the number of expansion coefficients of the Laguerre series for various values of the scaling parameter $\eta \in [50,1600]$. The calculations were made both with double real precision (Fig.~\ref{pic:test1}a,b) and single precision (Fig.\ref{pic:test1}c,d). To exclude the fictitious periodicity in algorithm 1, the initial approximation interval was increased to $t \in [0,2]$ , where $f(t)\equiv 0$ for $t \in [1,2]$. One can see in Fig.~\ref{pic:test1}a that an error of the order $\epsilon=10^{-14}$  was obtained with algorithm 1 for parameters $\eta=1600$ and $n=380\div920$, as well as for $\eta=800$ and $n=420\div440$. As the number of expansion coefficients for $n>920$ and $\eta=1600$ and for $n>440$ and $\eta=800$ increases, the approximation accuracy abruptly decreases, due to the fictitious periodicity and an abrupt break in the values of the series coefficients. As shown in Fig.~\ref{pic:test11}, the smaller is a given value of the parameter $\eta$, the longer is the spectrum. In this case the spectra of two periods of the function intersect at smaller values of $n$ and, starting with some number $n>n_0(\eta)$, the Laguerre series does not converge to the function being approximated.

In contrast to algorithm 1, the use of algorithm 2 (Fig.~\ref{pic:test1}b,d) did not require any additional increase in the approximation interval. The accuracy level of algorithm 2 is the same both in double and single real arithmetic and is of the order of $\epsilon=10^{-7}$, not decreasing to $\epsilon=10^{-14}$ as for algorithm 1. This is explained by the fact that when using formula (\ref{reverse_proc}) the sequence $l_m(\eta t_0), \ m=0,1,2... $ is the expansion coefficients for the delta function $\delta(t_0)$, for which (as shown in Fig.~\ref{pic:conjg24}) the Laguerre spectrum is infinitely long and slowly attenuating. Therefore, the finite number of expansion terms is a source of an additional error. However, the behavior of the error for algorithm 2 is more regular, since no fictitious periodicity and no abrupt break of the spectrum are observed for this calculation method.

\begin{figure}[!htb]
\centering
\begin{subfigure}[b]{0.47\textwidth}
\includegraphics[width=\textwidth]{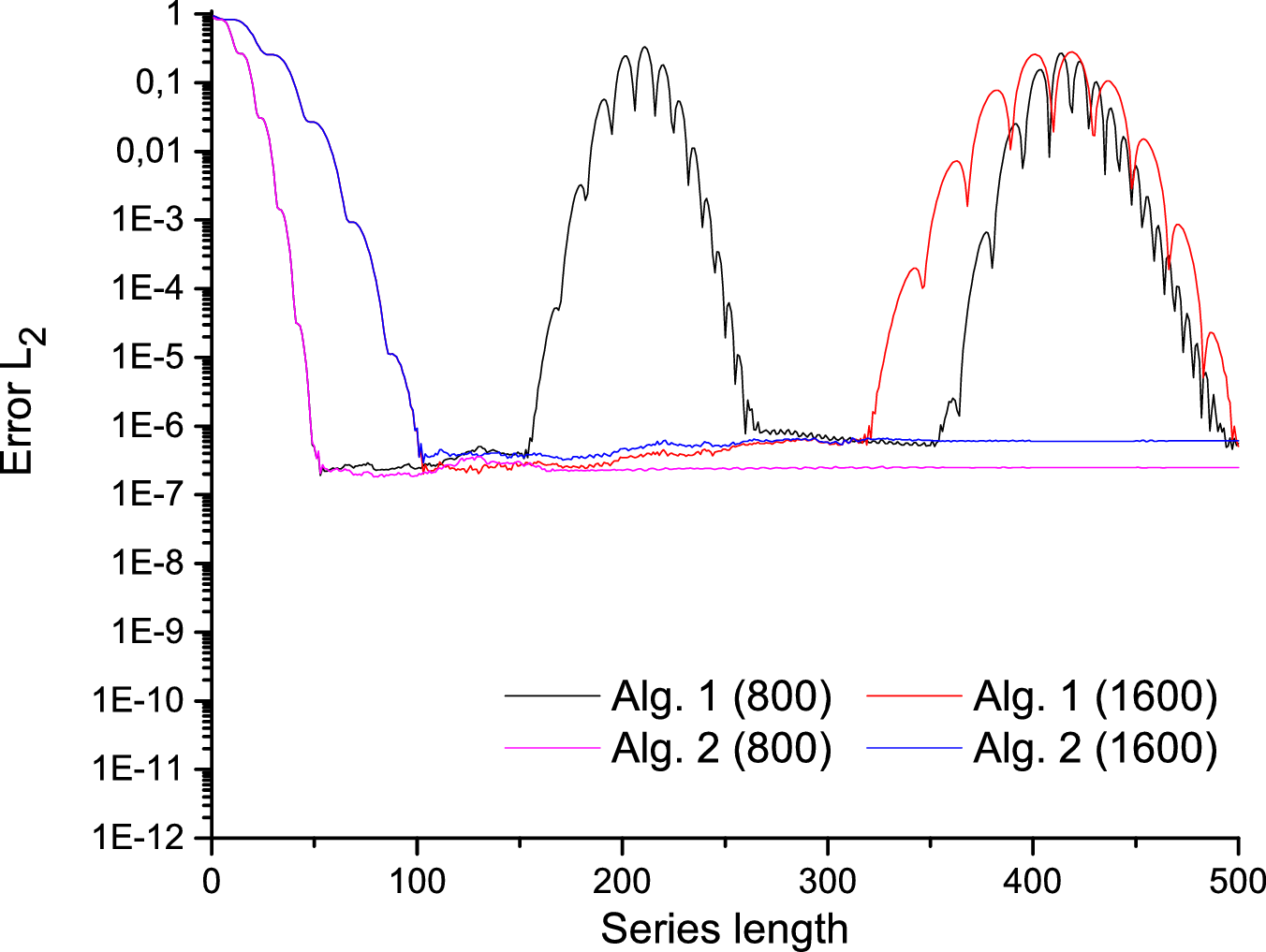}
\end{subfigure}
~
\begin{subfigure}[b]{0.47\textwidth}
\includegraphics[width=\textwidth]{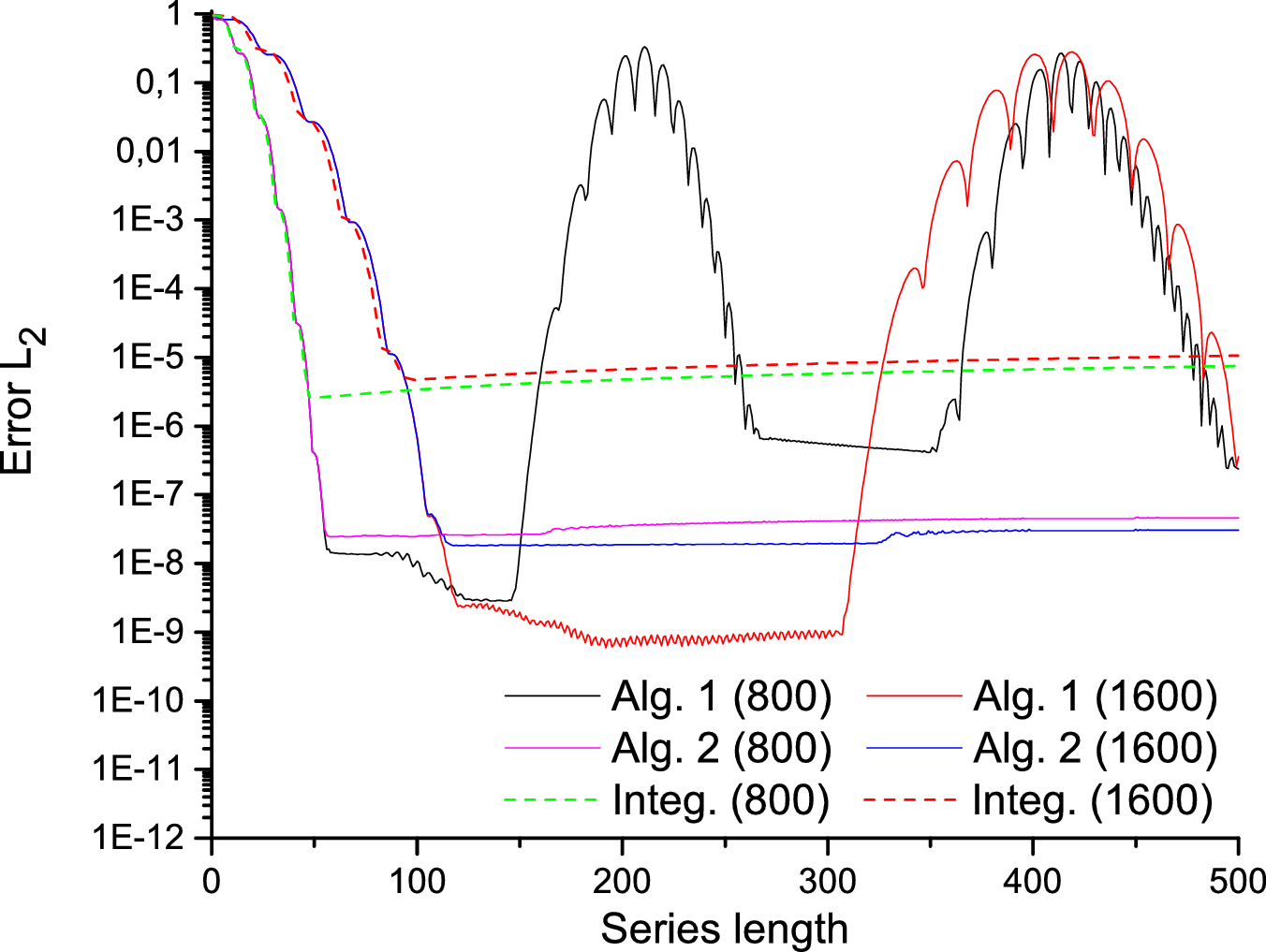}
\end{subfigure}
\caption{Approximation error versus the number of terms of Laguerre series for the function in Fig.~\ref{pic:conjg4}. Calculations were made by algorithms 1 and 2 for a) 32-bit precision, b) 64-bit precision.}
\label{pic:test12}
\end{figure}

Consider an approximation of the function in Fig.~\ref{pic:conjg4}, for which $f(0) = 0$. For this function, a comparison of the accuracy of algorithms 1 and 2 was made, with direct calculation of the integral (\ref{series_lag.int}) by the method of rectangles. Since the method of rectangles has the first order of accuracy, an integration step $4\times10^5$ times smaller than the discretization step for algorithms 1 and 2 was taken. Such a fine step is needed to provide high accuracy in calculating the Laguerre coefficients for a method of first order accuracy. To perform numerical integration of rapidly oscillating functions with a larger step, it is necessary to use quadratures of very high order of accuracy. In paper \cite{Litko1989} it was proposed to use quadratures of the $256$th order of accuracy to calculate a Laguerre series of length $n=128$. However, in solving practical problems nonsmooth functions have to be approximated. This calls into question whether it is reasonable to use high-accuracy quadratures for which error estimation implies the presence of high-order derivatives of the function to be expanded in the series.

By means of calculations with single real precision the calculation time can be decreased using a higher degree of vectorization of the calculations. In this case the error of algorithm 1 for $\eta=800,1600$ increases from $\epsilon=10^{-14}$ to $10^{-7}$, that is, to the level of single computer precision, and the error of algorithm 2 remains  at a level of the order of $\epsilon=10^{-7}$. It is important that the stability of all algorithms proposed still holds. Note that computer precision in calculations for nonsmooth functions may not be achieved, which eliminates the need for double real precision. To demonstrate this, in the test below we consider an approximation of a time series from a set of test seismograms for a velocity model called "Sigsbee"{}\cite{Paffenholz}.
\subsection{Test 2. Expansion of a non-smooth function}
\begin{figure}[!htb]
\centering
\begin{subfigure}[b]{0.47\textwidth}
\includegraphics[width=\textwidth]{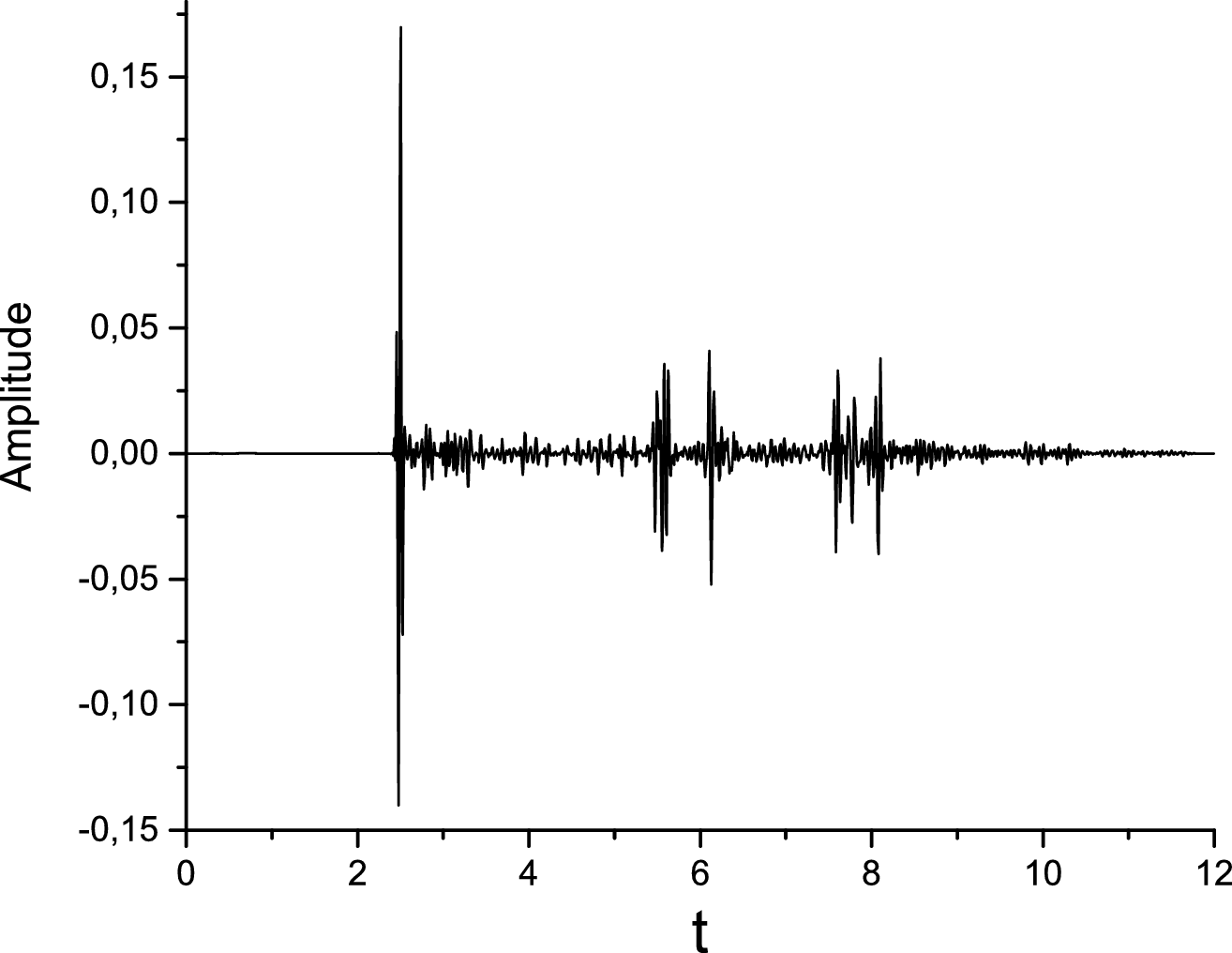}
\end{subfigure}
~
\begin{subfigure}[b]{0.47\textwidth}
\includegraphics[width=\textwidth]{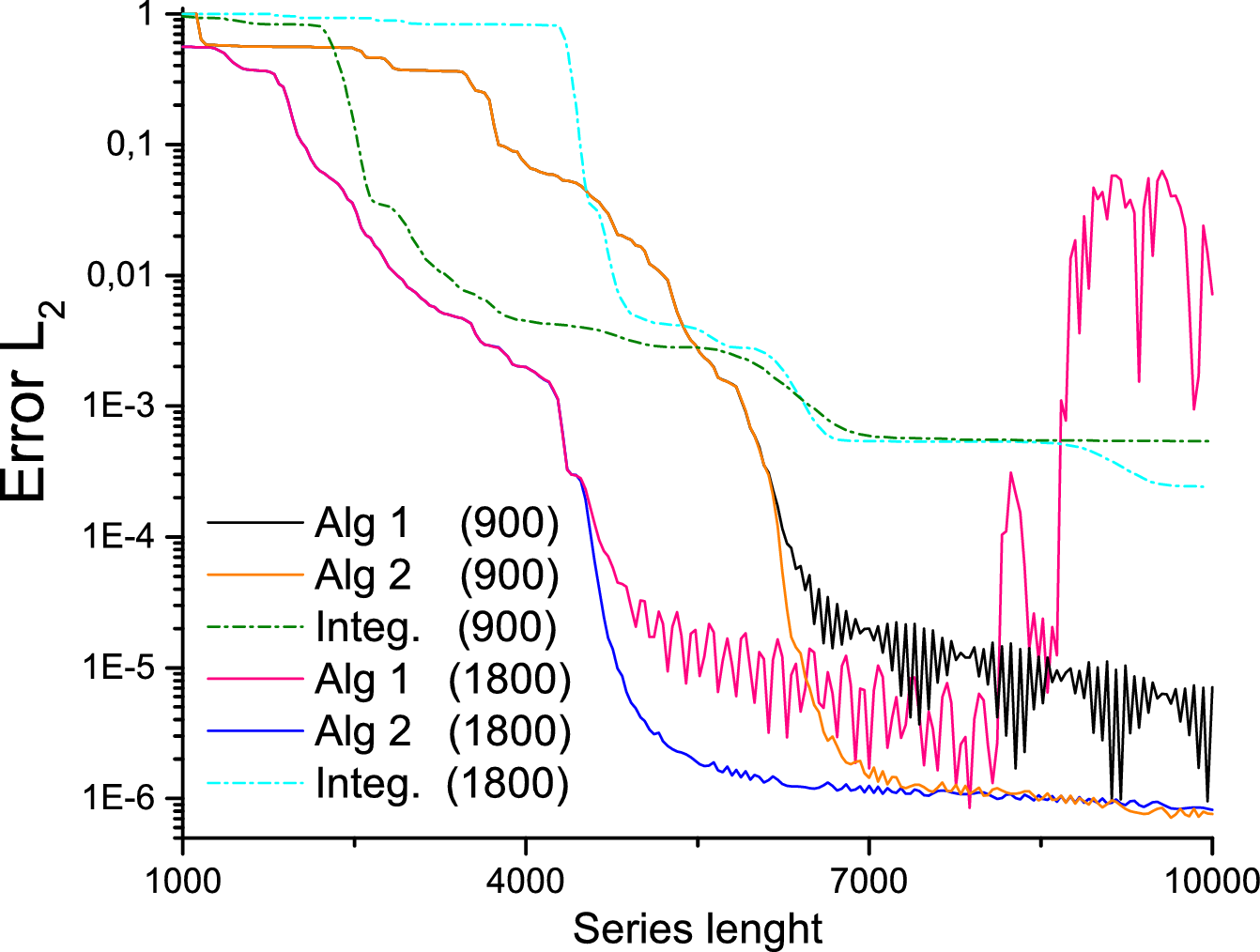}
\end{subfigure}
\caption{a) First trace from seismograms for the velocity model Sigsbee, b) approximation error versus the number of terms of the Laguerre series for the algorithms and transformation parameter $\eta=900,1800$ for first trace from seismograms for the velocity model Sigsbee.}
\label{pic:test2}
\end{figure}

In the first test an approximation of a smooth function on the interval $[0,1]$ was considered. Now let us test the above developed algorithms for a nonsmooth function (see Fig.~\ref{pic:test2}) specified on the interval $[0,12]$ with a discretization step $h_t=0.008$. This function corresponds to the first seismic trace from a test set of seismograms for the velocity model Sigsbee \cite{Paffenholz}. The seismograms for the SigSbee model have single real accuracy. Therefore, we consider an implementation of algorithms 1 and 2 only with single precision.

For calculations by algorithm 2 the approximation interval was not changed, whereas for algorithm 1 the approximation interval was increased by a factor of three up to $[0,36]$ by adding zero values (zero padding). One can see in Fig.~\ref{pic:test2} that algorithm 2 is approximately an order of magnitude more accurate than algorithm 1. Also, algorithm 2 demonstrates more regular behavior of the error, which considerably simplifies the process of finding an optimal number of coefficients of the Laguerre series. The smaller accuracy of algorithm 1 is caused, first, by the long Laguerre spectrum for the nonsmooth function, which does not make it possible to separate the spectra for different periods and exclude the influence of the fictitious periodicity. Second, the threefold increase in the approximation interval increases the number of terms of the series (\ref{main_formula}), which is a source of additional error due to the corresponding increase in the total number of operations.

The accuracy of calculating the expansion coefficients by the above proposed algorithms and that of calculating the integral (\ref{series_lag.int})  by the method of rectangles with an integration step $ht=8\times 10^{-7}$ was also compared. One can see in Fig.~\ref{pic:test2}b that, despite the fact that the discretization step of the function is much smaller, the accuracy of calculation by the method of rectangles is much lower than for algorithms 1 and 2, whereas the calculation burden is several orders of magnitude greater. It follows from formula (\ref{eq:assymptotic}) and Fig.~\ref{pic:laguerre_function} that as $n$ increases, the oscillation frequency of the functions $l_n(\eta t)$ also increases. Therefore, to calculate every subsequent expansion coefficient one has to either decrease the discretization step of the integrand or increase the order of the quadrature formula. However, when solving practical problems one should take into account that limited smoothness of the functions to be approximated may not allow using the maximum order of accuracy of the quadrature formula. Also, there are additional difficulties in using high-accuracy quadratures, which are caused by the need to calculate the integrand on a nonuniform grid, whereas a discrete function to be approximated is, as a rule, specified for equidistant values of the argument. In summary, we can say that the use of single precision for algorithms 1 and 2 is justified, since the observed error level is acceptable in solving practical problems. The use of double precision may decrease the error, but only in the approximation of very smooth functions.

\begin{table}[!h]
\center \small
\begin{tabular}{llllllllllllll}
  \hline
 & \multicolumn{6}{c}{$n=4096,\ t\in[0,12]$}& \multicolumn{6}{c}{$ n=8192,\ t\in[0,12]$}\\ \cmidrule(r){2-7} \cmidrule(r){8-13}
    p&$ {\mathrm{\epsilon}}$  & $\mathrm{Prec.}$&$ \mathrm{Step\ 1}$&$ \mathrm{Step\ 2}$&$ {\mathrm{Total}}$ & $\mathrm{Rel.}$&$ \mathrm{\epsilon}$&$ {\mathrm{Prec.}}$  & $\mathrm{Step\ 1}$&$ \mathrm{Step\ 2}$&$ \mathrm{Total}$&$ \mathrm{Rel.}$

    \\ \hline
 1  &1.7E-3&1.8&29.0&-&29.0&-&2.5E-6&3.9&52.0& -&52.0&-\\
  2  &2.9E-3&0.3&15.6 &0.9&16.5 &1.7&1.5E-5&1.1&27.6&2.1&29.7&1.7 \\
  3  &5.5E-3&5.7E-2&8.3&3.6&11.9&2.4&9.4E-5&0.2&15.1&9& 24.1&2.1\\
  8  &7.8E-3&1.2E-2&5.3&7.8&13.1&2.2&1.2E-4&3.9E-2&8.7&17&25.7& 2.0\\
  16  &2.2E-2&8.9E-2&3.5&15.5&19.1&1.5&2.1E-4&9.1E-3&5.7&34&39.7& 1.3\\
  \hline
  \end{tabular}
\caption{Estimates of calculation time and accuracy of algorithm 4. The number of auxiliary intervals {\it p} versus: \mbox{($\epsilon$)} approximation accuracy; \mbox{({\it Prec.})} preparatory calculation time needed to calculate local matrix $\tilde{M}$; \mbox{({\it Step 1})} calculation time  of local approximation; \mbox{({\it Step 2})} calculation time of the sequence of shifts for constructing the global approximation; \mbox{({\it Total})} total calculation time of algorithm 4; \mbox{({\it Rel.})} ratio between calculation time for algorithm 2 and calculation time for algorithm 4.}
\label{table33}
\end{table}

\begin{table}[!h]
\center \small
\begin{tabular}{llllllllllllll}
  \hline
 & \multicolumn{6}{c}{$n=32768,\ t\in[0,60]$}& \multicolumn{6}{c}{$ n=65536,\ t\in[0,120]$}\\ \cmidrule(r){2-7} \cmidrule(r){8-13}
    p&$ {\mathrm{\epsilon}}$  & $\mathrm{Prec.}$&$ \mathrm{Step\ 1}$&$ \mathrm{Step\ 2}$&$ {\mathrm{Total}}$ & $\mathrm{Rel.}$&$ \mathrm{\epsilon}$&$ {\mathrm{Prec.}}$  & $\mathrm{Step\ 1}$&$ \mathrm{Step\ 2}$&$ \mathrm{Total}$&$ \mathrm{Rel.}$
    \\ \hline
 1  &  5.3E-6& 77&877&-&877&-&               3.4E-6&359&3455& -&3455&-\\
  2  &8.1E-6&  19&449 &25&474&1.8&        4.6E-6&78.7&1754&58&1792&1.9 \\
  3  &6.1E-5&  4.8&231&44&275&3.1&        9.7E-6&19.4&826&107& 933&3.7\\
  8  &6.5E-5&  1.2&124&63&187&4.6&        6.6E-5&4.8&463&144&607& 5.7\\
  16  &5.8E-5& 2.5&79&80&159&5.5&         6.8E-5&1.2&242&180&422& 8.2\\
  32  &1.9E-4&  4.2E-2&52&91&143&6.1&   6.9E-5&2.6E-1&158&212&370& 9.3\\
  64  &5.8E-4&  9.6E-3&25&100&125&7.0& 1.2E-4&4.2E-2&105&239&344& 10\\
  128  &3.9E-3&  1.7E-3&17&112&129&6.8&5.8E-4&9.6E-3&50&265&315& 11\\
  \hline
  \end{tabular}
\caption{Estimates of calculation time and accuracy of algorithm 4. The number of auxiliary intervals {\it p} versus: \mbox{($\epsilon$)} approximation accuracy; \mbox{({\it Prec.})} preparatory calculation time needed to calculate local matrix $\tilde{M}$; \mbox{({\it Step 1})} calculation time  of local approximation; \mbox{({\it Step 2})} calculation time of the sequence of shifts for constructing the global approximation; \mbox{({\it Total})} total calculation time of algorithm 4; \mbox{({\it Rel.})} ratio between calculation time for algorithm 2 and calculation time for algorithm 4.}
\label{table34}
\end{table}

In testing of algorithm 4, Tables~\ref{table33} and \ref{table34} present the calculation times and accuracy estimates in the approximation of all $152684$ seismic traces for the Sigsbee model. For the interval $[0,12]$ the numbers of coefficients of the series were $n=4096$ and $8192$, and for the intervals $[0,60]$ and $[0,120]$ the numbers of coefficients of the series were specified as $n=32768$ and $65536$, respectively. Initial seismic traces for the Sigsbee model were specified for $t \in [0,12]$. To obtain the time series for the intervals $ [0,60]$ and $[0,120]$, the initial trace was supplemented by four or eleven identical copies of the initial signal, respectively. One can see from the data presented that, although algorithm 4 does not belong to the class of fast algorithms, it allows a slight decrease in the calculation time, especially for large time intervals. Also note that when using algorithm 4 the time of the preparatory calculations needed to modify the matrix $M$ in implementing algorithm 2 decreases considerably, since a matrix of smaller order is required to approximate the local functions. Nevertheless, it follows from Tables~\ref{table33} and \ref{table34} that the approximation accuracy $\epsilon$ decreases as the number of auxiliary intervals increases. This is caused, first, by the presence of auxiliary buffers, in which multiplication by an exponentially decreasing factor is made for the Laguerre spectrum of a local function not to be infinite because of the discontinuities of the function values on the boundaries of the subintervals. On each subinterval the local function is approximated by a Laguerre series with a number of coefficients of $n/p$, where $p$ is the number of subintervals. However, $n/p$ expansion coefficients may be insufficient to approximate a nonsmooth local function, which results in loss in approximation accuracy. Nevertheless, if one has to approximate a time series with an accuracy of the order of $\epsilon=10^{-3}\div10^{-5}$ (which is sufficient for practical calculations  \cite{Terekhov2017,Terekhov2018}), it is recommended to use algorithm 4.

To multiply the matrix $\tilde{M}$ by a vector, a numerical procedure from BLAS MKL library was used. Taking into account high degree of optimization of the BLAS procedure for a specific processor model, the calculation of Laguerre coefficients is performed very fast. At the same time, optimization of the algorithm of fast Fourier transform is a more complicated problem. This decreases the degree of vectorization of the calculations at the second step of algorithm 4. As a result, the speedup of algorithm 4 also decreases. If an internal Fortran procedure, such as "matlmul"{}, had been used for matrix multiplication instead of that from BLAS library, the speedup coefficient of algorithm 4 would have been much larger (although the total calculation time also increases), since the computational costs of the "matmul"{} function are several times greater than those of the procedure from BLAS MKL library.
\section{Conclusions}
In this paper, new algorithms to calculate the integral Laguerre transform by solving a one-dimensional transport equation have been developed. The main idea of the above proposed approach is that the calculation of improper integrals of rapidly oscillating functions is replaced by solving an initial boundary value problem for the transport equation using spectral algorithms. This approach has made it possible to successfully avoid the problems formulated in the introduction and associated with numerical implementation of the Laguerre transform. It would have been impossible to implement the above proposed computational model without the development of auxiliary procedures that allow removing the fictitious periodicity resulting from periodic boundary conditions. One of the correcting procedures is based on solving the transport equation, whereas the other one is based on a posteriori analysis of the Laguerre spectrum energy. Test calculations have shown that the first method of removing the periodicity is more reliable, accurate, and efficient, since it does not require increasing the approximation interval. Although the above algorithms do not belong to the class of fast algorithms, the number of arithmetic operations has been considerably decreased, since there is no need in calculations with small grid steps or quadrature formulas of high orders of accuracy to calculate rapidly oscillating improper integrals. Additionally, an approach has been developed to decrease the computational costs in making the Laguerre transform for large approximation intervals by solving the transport equation. The test calculations have also confirmed that all developed algorithms can be used both with single and double real precision without loss of numerical stability. Thus, if a large set of functions is approximated by a Laguerre series (for instance, in solving problems of seismic prospecting), the above proposed algorithms allow saving the calculation time considerably. This fact makes this approach attractive from both theoretical and practical viewpoints.
\section{Acknowledgements}
The work was partially supported by a grant from the Russian Foundation for Basic Research RFBR, grant \mbox {no. 18-41-543002}.
\newpage
\bibliography{base}

\begin{thebibliography}{10}

\bibitem{fatab2011}
A.~G. Fatyanov and A.~V. Terekhov.
\newblock High-performance modeling acoustic and elastic waves using the
  parallel dichotomy algorithm.
\newblock {\em J. Comp. Phys.}, 230(5):1992--2003, 2011.

\bibitem{Terekhov:2013}
A.~V. Terekhov.
\newblock A fast parallel algorithm for solving block-tridiagonal systems of
  linear equations including the domain decomposition method.
\newblock {\em Parallel Comput.}, 39(6-7):245--258, 2013.

\bibitem{Terekhov2015206}
A.~V. Terekhov.
\newblock Spectral-difference parallel algorithm for the seismic forward
  modeling in the presence of complex topography.
\newblock {\em Journal of Applied Geophysics}, 115(0):206--219, 2015.

\bibitem{Mikhailenko1999}
B.~G. Mikhailenko.
\newblock Spectral {Laguerre} method for the approximate solution of time
  dependent problems.
\newblock {\em Applied Mathematics Letters}, 12:105--110, 1999.

\bibitem{Mikhailenko2008}
B.~G. Mikhailenko and A.~F. Mastryukov.
\newblock Numerical solution of maxwell's equations for anisotropic media using
  the laguerre transform.
\newblock {\em Russian Geology and Geophysics}, 49:621--627, 2008.

\bibitem{Colton1984}
D.~Colton and J.~Wimp.
\newblock Analytic solutions of the heat equation and some formulas for
  {L}aguerre and {H}ermite polynomials.
\newblock {\em Complex Variables, Theory and Application: An International
  Journal}, 3(4):397--412, 1984.

\bibitem{Jo2006}
Javier A.~{et al} Jo.
\newblock Laguerre-based method for analysis of time-resolved fluorescence
  data: Application to in-vivo characterization and diagnosis of
  atherosclerotic lesions.
\newblock {\em Journal of biomedical optics}, 11(2), 2006.

\bibitem{Terekhov2017}
A.~V. Terekhov.
\newblock The {Laguerre} finite difference one-way equation solver.
\newblock {\em Computer Physics Communications}, 214:71 -- 82, 2017.

\bibitem{Terekhov2018}
A.~V. Terekhov.
\newblock The stabilization of high-order multistep schemes for the {Laguerre}
  one-way wave equation solver.
\newblock {\em Journal of Computational Physics}, 368:115 -- 130, 2018.

\bibitem{Weeks1966}
W.~T. Weeks.
\newblock Numerical inversion of {L}aplace transforms using {L}aguerre
  functions.
\newblock {\em J. ACM}, 13(3):419--429, July 1966.

\bibitem{Abate1996}
J.~Abate, G.~Choudhury, and W.~Whitt.
\newblock On the laguerre method for numerically inverting laplace transforms.
\newblock {\em Informs J. on Computing}, 8(4):413--427, 1996.

\bibitem{Strain1992}
J.~Strain.
\newblock A fast {Laplace} transform based on {Laguerre} functions.
\newblock {\em Mathematics of Computation}, 58(197):275--283, 1992.

\bibitem{Weber1980}
H.~Weber.
\newblock Numerical computation of the fourier transform using laguerre
  functions and the fast {Fourier} transform.
\newblock {\em Numerische Mathematik}, 36(2):197--209, Jun 1980.

\bibitem{Golub1989}
G.~H. Golub and C.~F. Van~Loan.
\newblock {\em Matrix computations (3rd ed.)}.
\newblock Johns Hopkins University Press, Baltimore, MD, USA, 1996.

\bibitem{Samarski_Nikolaev}
A.A. Samarskij and E.S. Nikalayev.
\newblock {\em Numerical Methods for Grid Equations}.
\newblock Birkhauser Verlag, 1989.

\bibitem{terekhov:Dichotomy}
A.~V. Terekhov.
\newblock Parallel dichotomy algorithm for solving tridiagonal system of linear
  equations with multiple right-hand sides.
\newblock {\em Parallel Comput.}, 36(8):423--438, 2010.

\bibitem{Terekhov2016}
A.~V. Terekhov.
\newblock A highly scalable parallel algorithm for solving toeplitz tridiagonal
  systems of linear equations.
\newblock {\em Journal of Parallel and Distributed Computing}, 87:102--108,
  2016.

\bibitem{abramowitz+stegun}
M.~Abramowitz and I.~A. Stegun.
\newblock {\em Handbook of Mathematical Functions with Formulas, Graphs, and
  Mathematical Tables}.
\newblock Dover, New York, ninth dover printing, tenth gpo printing edition,
  1964.

\bibitem{Szegoe1975}
G.~Szeg{\"o}.
\newblock {\em Orthogonal Polynomials}.
\newblock American Math. Soc: Colloquium publ. American Mathematical Society,
  1975.

\bibitem{Rainville1971}
D.~Rainville.
\newblock {\em Special Functions}.
\newblock AMS Chelsea Publishing Series. Chelsea Publishing Company, 1971.

\bibitem{Temme1990}
N.~M. Temme.
\newblock Asymptotic estimates for {Laguerre} polynomials.
\newblock {\em Zeitschrift f{\"u}r angewandte Mathematik und Physik ZAMP},
  41(1):114--126, Jan 1990.

\bibitem{Litko1989}
Joseph~R. Litko.
\newblock Gi/g/1 interdeparture time and queue-length distributions via the
  {Laguerre} transform.
\newblock {\em Queueing Systems}, 4(4):367--381, Dec 1989.

\bibitem{Nussbaumer1982}
H.~J. Nussbaumer.
\newblock {\em Fast Fourier Transform and Convolution Algorithms}.
\newblock Springer-Verlag, 1982.

\bibitem{Boyd2001}
J.~P. Boyd.
\newblock {\em Chebyshev and Fourier Spectral Methods}.
\newblock Dover, New York, 2001.

\bibitem{GOHBERG1994411}
I.~Gohberg and V.~Olshevsky.
\newblock Fast algorithms with preprocessing for matrix-vector multiplication
  problems.
\newblock {\em Journal of Complexity}, 10(4):411 -- 427, 1994.

\bibitem{Pan2016}
V.~Pan.
\newblock How bad are vandermonde matrices?
\newblock {\em SIAM Journal on Matrix Analysis and Applications},
  37(2):676--694, 2016.

\bibitem{Gautschi2011}
W.~Gautschi.
\newblock Optimally scaled and optimally conditioned {V}andermonde and
  {V}andermonde-like matrices.
\newblock {\em BIT Numerical Mathematics}, 51(1):103--125, Mar 2011.

\bibitem{Borodin1974}
A.~Borodin and R.~Moenck.
\newblock Fast modular transforms.
\newblock {\em Journal of Computer and System Sciences}, 8(3):366 -- 386, 1974.

\bibitem{Gohberg1994}
I.~Gohberg and V.~Olshevsky.
\newblock Complexity of multiplication with vectors for structured matrices.
\newblock {\em Linear Algebra and its Applications}, 202:163 -- 192, 1994.

\bibitem{Pan1993}
V.~Pan, A.~Sadikou, E.~Landowne, and O.~Tiga.
\newblock A new approach to fast polynomial interpolation and multipoint
  evaluation.
\newblock {\em Computers \& Mathematics with Applications}, 25(9):25 -- 30,
  1993.

\bibitem{Gautschi1983}
W.~Gautschi.
\newblock The condition of vandermonde-like matrices involving orthogonal
  polynomials.
\newblock {\em Linear Algebra and its Applications}, 52-53:293 -- 300, 1983.

\bibitem{Alpert1991}
B.~Alpert and V.~Rokhlin.
\newblock A fast algorithm for the evaluation of {Legendre} expansions.
\newblock {\em SIAM Journal on Scientific and Statistical Computing},
  12(1):158--179, 1991.

\bibitem{Hale2016}
N.~Hale and A.~Townsend.
\newblock A fast {FFT-based} discrete {Legendre} transform.
\newblock {\em IMA Journal of Numerical Analysis}, 36(4):1670--1684, 2016.

\bibitem{Leibon2008}
G.~Leibon, D.~N. Rockmore, W.~Park, R.~Taintor, and G.~S. Chirikjian.
\newblock A fast {Hermite} transform.
\newblock {\em Theoretical Computer Science}, 409(2):211 -- 228, 2008.
\newblock Symbolic-Numerical Computations.

\bibitem{Bostan2010}
A.~Bostan, B.~Salvy, and E.~Schost.
\newblock Fast conversion algorithms for orthogonal polynomials.
\newblock {\em Linear Algebra and its Applications}, 432(1):249 -- 258, 2010.

\bibitem{ONeil2010}
M.~O'Neil, F.~Woolfe, and V.~Rokhlin.
\newblock An algorithm for the rapid evaluation of special function transforms.
\newblock {\em Applied and Computational Harmonic Analysis}, 28(2):203 -- 226,
  2010.
\newblock Special Issue on Continuous Wavelet Transform in Memory of Jean
  Morlet, Part I.

\bibitem{Smith1985}
D.~R. Smith.
\newblock The design of divide and conquer algorithms.
\newblock {\em Science of Computer Programming}, 5:37 -- 58, 1985.

\bibitem{Gil2017}
A.~Gil, J.~Segura, and N.~M. Temme.
\newblock Efficient computation of laguerre polynomials.
\newblock {\em Computer Physics Communications}, 210:124 -- 131, 2017.

\bibitem{Paffenholz}
J.~Paffenholz, B.~McLain, J.~Zaske, and P.~J. Keliher.
\newblock {\em Subsalt multiple attenuation and imaging: Observations from the
  Sigsbee2B synthetic dataset}, chapter 538, pages 2122--2125.

\end{thebibliography}
\end{document}